\documentclass[11pt]{article}

\usepackage[margin=1in]{geometry}
\usepackage{amsmath,amssymb,amsthm,mathtools,bm}
\usepackage{bbm}
\usepackage{booktabs}
\usepackage{graphicx}
\usepackage{enumitem}
\usepackage{hyperref}
\usepackage[round,authoryear]{natbib}
\usepackage{xcolor}
\usepackage{array}
\usepackage{caption}
\usepackage{subcaption}
\usepackage{microtype}
\usepackage{algorithm,algpseudocode}

\hypersetup{
    colorlinks=true,
    linkcolor=blue!60!black,
    citecolor=blue!60!black,
    urlcolor=blue!60!black
}

\newtheorem{theorem}{Theorem}[section]
\newtheorem{proposition}[theorem]{Proposition}
\newtheorem{lemma}[theorem]{Lemma}
\newtheorem{corollary}[theorem]{Corollary}
\theoremstyle{definition}

\newtheorem{remark}[theorem]{Remark}

\newcommand{\Unif}{\mathrm{Unif}}

\title{Reference-Null Calibrated Thresholds for E-Processes with Applications to Conformal Martingales}
\author{
Yi Ding$^{1}$,
Lan Wei$^{2}$
Xuehu Zhu$^{1,*}$,
and Wenlin Dai$^{2,*}$
}
\date{}
\begin{document}
\maketitle

\begin{center}
\small
$^{1}$ School of Mathematics and Statistics, Xi'an Jiaotong University, Xi'an, China\\
$^{2}$ Institute of Statistics and Big Data, Renmin University of China, Beijing, China
\end{center}
\begingroup
\renewcommand{\thefootnote}{}
\footnotetext{\textsuperscript{*} Co-corresponding authors.}
\addtocounter{footnote}{-1}
\endgroup

\begin{abstract}
    E-processes provide a flexible framework for anytime-valid inference, 
    with the conventional rejection boundary $1/\alpha$ typically justified by 
    Ville's inequality. Such a boundary is universal but can be conservative, 
    as it does not exploit additional information about the null distribution. 
    We propose a reference-null calibration framework that uses independent null 
    samples to construct sharper rejection thresholds while preserving type-I 
    error control. 
To study the statistical gain from sharper thresholds, we specialize the
framework to conformal martingales. Building on histogram-based betting, we
incorporate Krichevsky--Trofimov smoothing and develop a restart-mixture
construction for distribution shift detection. We establish
quantitative results for detection power and detection delay and characterize
how the reduction in the rejection boundary translates into power and delay
gains relative to the conventional Ville boundary.
 Numerical experiments corroborate the theoretical 
findings across a range of distributional changes. Finally, we apply the 
proposed calibration strategy to an existing e-process for online LLM watermark 
detection, showing that the method can improve sequential detection performance 
without modifying the underlying e-process. These results demonstrate that 
reference-null calibration provides a general and modular way to enhance the 
efficiency of e-process-based sequential inference.

\end{abstract}

\noindent\textbf{Keywords:} e-values; e-processes; conformal martingales; Ville's inequality; Markov's inequality; anytime-valid inference; reference-null calibration; p-to-e calibrators.

\section{Introduction}

Sequential statistical inference concerns hypothesis testing and estimation
problems in which observations arrive sequentially and statistical evidence is
continuously updated over time. Classical statistical procedures are typically
developed for a fixed sample size and rely on a pre-specified stopping rule.
However, many modern applications require monitoring data streams and making
decisions as soon as sufficient evidence has accumulated, without knowing in
advance when data collection will stop. Examples include online experiments,
clinical monitoring, and sequential quality control, where early decisions can
reduce cost and improve resource allocation
\citep{wald1945sequential,lai1988nearly,kharitonov2015sequential,lindon2024anytime}.
Such settings create a fundamental statistical challenge: repeatedly examining
accumulating data can inflate type-I error unless the inference procedure
remains valid under arbitrary stopping rules. This challenge has motivated the
development of anytime-valid inference, which provides error guarantees that
hold uniformly over time and under optional stopping
\citep{ramdas2023game,waudby2023distribution}.

A recent framework addressing this challenge is provided by e-values and
e-processes, which have emerged as a flexible foundation for anytime-valid
inference. The idea of e-values originates from the game-theoretic approach to
statistics and provides a nonnegative measure of evidence whose expectation is
controlled under the null hypothesis
\citep{shafer2021testing,vovk2021evalues}. Specifically, an e-variable
satisfies a unit expectation constraint under the null, allowing Markov's
inequality to provide a valid rejection rule without requiring a fixed sample
size. E-processes extend this principle to sequential settings by requiring
validity at arbitrary stopping times. They therefore provide a unified
framework for optional stopping, continuous monitoring, and adaptive data
collection; see \citet{ramdas2023game,grunwald2024safe} for recent
developments and perspectives on safe anytime-valid inference.

The universal validity of e-processes, however, comes with a potential loss of
statistical efficiency. For an e-process $(M_t)_{t\geq0}$, Ville's inequality
gives the anytime-valid guarantee
\[
\mathbb{P}_0\left(
\sup_{t\geq0}M_t\geq \frac1\alpha
\right)\leq\alpha .
\]
The resulting boundary $1/\alpha$ is attractive because it is distribution-free
and applies to every nonnegative e-process. Nevertheless, it is derived from a
worst-case argument and does not exploit the specific null distribution of the
particular e-process under consideration. In many applications, the actual
crossing probability of a given e-process is substantially smaller than the
upper bound provided by Ville's inequality, leaving part of the available
type-I error budget unused. This conservativeness motivates the development of
sharper rejection boundaries that preserve anytime validity while improving
statistical power.

Several recent studies have investigated ways to improve the efficiency of
e-values and e-processes. One line of work derives sharper thresholds by
imposing additional structural assumptions on the null distribution of the
e-variable, such as shape constraints or restrictions on transformed
densities
\citep{blierwong2024improved}. These methods demonstrate that the Markov
threshold can be substantially improved when the null distribution belongs to
a suitable restricted class. Another direction focuses on refining the
behavior of specific sequential processes. For example,
\citet{fischer2026improving} study improvements based on boundary overshoot in
likelihood-ratio-type e-processes by modifying the evidence accumulation
mechanism, while \citet{delapena2026exact} characterize the slack in Ville's
inequality through an exact decomposition of boundary-crossing probabilities.
Although these approaches provide important insights, they either require
structural knowledge of the null process or modify the original e-process
itself. A general procedure for improving the rejection boundary of an
already constructed e-process, using only external information about its null
behavior, remains desirable.

In this paper, we introduce a reference-null calibration framework for
constructing sharper rejection boundaries for e-values and e-processes. The
central idea is simple: instead of relying on a universal worst-case boundary,
we use independent reference-null realizations to directly estimate the null
distribution of the relevant evidence statistic. For a finite-horizon
e-process $(M_t)_{0\leq t\leq T}$, the relevant statistic is the path maximum
\[
S_T=\max_{0\leq t\leq T}M_t,
\]
because rejection occurs precisely when this quantity exceeds a prescribed
boundary. By generating independent null trajectories of the same e-process,
we estimate the upper quantile of $S_T$ and use it as a calibrated rejection
threshold. This approach has three advantages. First, it preserves the
original e-process and therefore separates boundary calibration from evidence
construction. Second, it does not require an analytically tractable null
distribution of the crossing statistic. Third, it provides finite-sample
validity through exchangeability between the test trajectory and the
reference-null trajectories. The use of exchangeable reference samples also connects our approach to the
general calibration principle underlying conformal inference
\citep{vovk2005algorithmic}.
However, our goal is different: instead of calibrating predictive sets or
p-values, we calibrate the rejection boundary of an already constructed
e-process.

To demonstrate the practical value of the proposed framework, we specialize
the calibration procedure to conformal martingales, a particularly natural
class of e-processes for sequential distribution-shift detection. Conformal
martingales were introduced as distribution-free sequential tests of
exchangeability
\citep{vovk2003testing,vovk2005algorithmic} and have since been extended through
adaptive betting strategies and plug-in constructions
\citep{fedorova2012plugin,eliades2020histogram,shaer2026cctm}. Their distribution-free null
behavior makes them an ideal setting for studying the effect of boundary
calibration independently of model misspecification. Within this application,
we consider histogram-based betting with Krichevsky--Trofimov smoothing and
develop a restart-mixture construction to address unknown distributional
change points. These components allow us to investigate how reference-null
calibration interacts with adaptive evidence accumulation and sequential
change detection.

The main contributions of this paper are summarized as follows.

\begin{enumerate}

\item We introduce a general reference-null calibration framework for improving
rejection boundaries of e-values and e-processes. Unlike approaches that modify
the underlying evidence process or require analytical characterization of its
null distribution, our method constructs sharper boundaries from independent
reference-null trajectories while preserving the validity of the original
e-process. This provides a modular approach for improving existing
anytime-valid inference procedures.

\item We establish theoretical guarantees for the efficiency gains achieved by
reference-null calibration. In particular, by applying the framework to
conformal martingales, we provide a systematic finite-sample analysis of
detection power and detection delay for adaptive histogram-based procedures
with and without restart aggregation. Existing studies of conformal martingales have largely emphasized validity
guarantees and empirical evaluations, while systematic finite-sample analyses
of detection power and delay remain limited. Our results characterize how calibrated boundaries
translate into sequential detection gains and quantify the effects of signal
strength, adaptation cost, and restart aggregation.

\item We demonstrate the practical value of reference-null calibration through
numerical experiments and an online LLM watermark detection application. The
application highlights the modularity of the proposed framework: an existing
e-process can be improved by recalibrating only its rejection boundary, without
altering the underlying evidence construction mechanism.

\end{enumerate}

The remainder of the paper is organized as follows. Section~\ref{sec:reference-null-calibration}
develops the reference-null calibration framework, beginning with e-values and
then extending the construction to finite-horizon e-processes. Section
\ref{sec:conformal-martingale} introduces the conformal martingale procedures
based on histogram betting and restart aggregation. Sections~4 and~5 establish
power and detection-delay guarantees. Section~6 presents numerical experiments
and an application to online LLM watermark detection. Section~7 concludes with
a discussion of limitations and future directions.
\section{Reference-Null Calibration for E-Values and E-Processes}
\label{sec:reference-null-calibration}

\subsection{Reference-null calibration for e-values: illustrative results}
\label{subsec:evalue-warmup}

We begin with the static e-value setting, which isolates the basic idea of
reference-null calibration without the additional complication of sequential
monitoring. Throughout this paper, we focus on simple null hypotheses, or on
problems that can be reduced to a simple null through an appropriate
transformation, conditioning argument, or pivotal representation. Accordingly,
let \(P_0\) denote the null distribution and let \(E\) be a nonnegative
e-variable satisfying
$
    \mathbb E_{P_0}[E]
    \leq
    1.
$

Let
$
    F_0(x)=P_0(E\leq x)
$
denote the null distribution function of \(E\). If \(F_0\) were known, the
natural distribution-aware rejection threshold would be its upper
\(\alpha\)-quantile,
\begin{equation}
    \label{eq:oracle-evalue-threshold}
    c_{\alpha}^{\star}
    :=
    F_0^{\leftarrow}(1-\alpha)
    :=
    \inf
    \left\{
        x\geq0:
        F_0(x)\geq1-\alpha
    \right\}.
\end{equation}
By the definition, we have
$
    P_0
    \left(
        E>c_{\alpha}^{\star}
    \right)
    \leq
    \alpha.
$
By Markov's inequality, we also have
$
    P_0
    \left(
        E\geq \frac{1}{\alpha}
    \right)
    \leq
    \alpha.
$
By the definition of \(c_{\alpha}^{\star}\) again, we have
$
    c_{\alpha}^{\star}
    \leq
    \frac{1}{\alpha}.
$
Hence \(c_{\alpha}^{\star}\), rather than \(1/\alpha\), is a
smaller threshold.
Markov's inequality is sharp over the class of all nonnegative random
variables. Any systematic improvement
over
\(1/\alpha\) must therefore exploit information beyond the defining
e-variable constraint. Shape-constrained approaches use structural
assumptions on \(F_0\) for this purpose
\citep{blierwong2024improved}; here, instead, we use independent
reference-null observations.

In practice, \(F_0\) and hence \(c_{\alpha}^{\star}\) are typically
unknown. Suppose that
$
    E_1^{\mathrm{ref}},
    \ldots,
    E_B^{\mathrm{ref}}
$
are reference e-values generated under \(P_0\), and let
\(E^{\mathrm{test}}\) denote the e-value computed from the test observation.
We require$
    \left(
        E^{\mathrm{test}},
        E_1^{\mathrm{ref}},
        \ldots,
        E_B^{\mathrm{ref}}
    \right)
$to be exchangeable under $P_0$.

Write
$
    E_{(1)}^{\mathrm{ref}}
    \leq
    \cdots
    \leq
    E_{(B)}^{\mathrm{ref}}
$
for the order statistics of the reference e-values, and adopt the convention
$
    E_{(B+1)}^{\mathrm{ref}}
    :=
    +\infty.
$
For a prescribed level \(\alpha\in(0,1)\), define
$
    \label{eq:evalue-calibration-index}
    k_{\alpha,B}
    :=
    \left\lceil
        (1-\alpha)(B+1)
    \right\rceil
$
and the reference-null calibrated threshold
$
    \widehat c_{\alpha,B}
    :=
    E_{(k_{\alpha,B})}^{\mathrm{ref}}.
$
We reject the null whenever
\(\phi_{\alpha,B} := \mathbf 1 \left\{ E^{\mathrm{test}} > \widehat c_{\alpha,B} \right\} = 1.\)

\begin{theorem}[Finite-sample validity of reference-null calibration]
    \label{thm:evalue-reference-validity}
    Fix \(\alpha\in(0,1)\) and \(B\in\mathbb N\). Suppose that
    $E^{\mathrm{test}}$ is an e-variable for \(P_0\) and that
    $\left(
        E^{\mathrm{test}},
        E_1^{\mathrm{ref}},
        \ldots,
        E_B^{\mathrm{ref}}
    \right)$ are exchangeable. Then
    \(\mathbb E_{P_0} \left[ \phi_{\alpha,B} \right] \leq \alpha.\)
    If the common null distribution of
    \(E^{\mathrm{test}},E_1^{\mathrm{ref}},\ldots,E_B^{\mathrm{ref}}\)
    is continuous, then it also holds that
    \(\mathbb E_{P_0} \left[ \phi_{\alpha,B} \right] \geq \alpha-\frac{1}{B+1}.\)
\end{theorem}

Theorem~\ref{thm:evalue-reference-validity} is a finite-sample statement:
no asymptotic approximation and no regularity condition on the distribution
of \(E\) are required. The result is driven by the exchangeability between
the test e-value and the reference e-values, rather than by the moment
condition defining an e-variable. This rank-based calibration argument is
closely related to the basic principle underlying conformal prediction, where
exchangeability between a test observation and a reference sample is used to
obtain exact finite-sample validity through the rank of the test statistic
\citep{vovk2005algorithmic,lei2018distribution,barber2021predictive}. Our use of this principle is
different in purpose: rather than constructing a prediction set, we use the
exchangeable rank to calibrate the rejection threshold of an already
constructed e-value.

This e-value warm-up highlights the basic distinction underlying our
approach. Markov's inequality controls the upper tail of \(E\) indirectly
through the single moment restriction \(\mathbb E_{P_0}[E]\leq1\), whereas
reference-null calibration estimates the relevant upper tail directly.
For an e-process, the same principle applies after replacing the scalar
e-value \(E\) by the finite-horizon path maximum of the process.

\subsection{Reference-null calibration of finite-horizon e-process boundaries}
\label{subsec:eprocess-calibration}

We next generalize the reference-null calibration framework from e-values to finite-horizon e-processes. Although Ville's inequality provides a universal
boundary for any nonnegative e-process, it does not exploit the particular
null distribution of the process under consideration. The central idea of
this section is to replace the worst-case Ville boundary by a
finite-horizon boundary calibrated from independent null trajectories.

Let
$
    (M_t)_{0\leq t\leq T}
$
be an e-process with respect to a filtration
\((\mathcal F_t)_{t\geq0}\) under the null distribution \(P_0\).
For a finite monitoring horizon \(T\), define the path maximum
$
    S_T
    :=
    \max_{0\leq t\leq T}M_t .
$
The usual Ville boundary rejects the null whenever
$
    S_T\geq \frac1\alpha ,
$
which guarantees
$
    P_0
    \left(
        S_T\geq\frac1\alpha
    \right)
    \leq\alpha .
$
The validity of this boundary follows solely from the e-process property
and therefore holds uniformly over all possible null distributions.
However, this universality also implies that the boundary is calibrated
against the least favorable class of nonnegative processes rather than
against the actual null behavior of the given e-process.

Suppose that independent reference-null trajectories are available. Such
trajectories may arise from historical null data, independent control
samples, or a data-generating mechanism under the null hypothesis. The
purpose of these trajectories is not to construct the e-process itself,
but only to estimate the finite-horizon crossing distribution of the
already specified e-process.

More generally, let
\(g=(g_0,\ldots,g_T)\)
be a deterministic positive boundary template, specified independently of
the calibration trajectories. The template allows one to incorporate
known temporal features of the monitoring problem. For example, one may
choose \(g_t\equiv1\), which corresponds to calibrating a constant
boundary for the path maximum, or use a nonconstant template when
different monitoring times require different scaling.

For any trajectory define the normalized path maximum
$
    R_g
    :=
    \max_{0\leq t\leq T}
    \frac{M_t}{g_t}.
$
Let
$
    R_{g,(1)},\ldots,R_{g,(B)}
$
be the corresponding statistics computed from \(B\) independent
reference-null trajectories. Denote their order statistics by
\(R_{g,(1)} \leq\cdots\leq R_{g,(B)} ,\)
and adopt the convention $R_{g,(B+1)}:=+\infty.$
For a significance level \(\alpha\), define
$
    k_{\alpha,B}
    :=
    \left\lceil
        (1-\alpha)(B+1)
    \right\rceil ,
$
and the calibrated multiplier
$
    \widehat q_{\alpha}(g)
    :=
    R_{g,(k_{\alpha,B})}.
$
The resulting rejection rule is
\[
    \nu_\alpha
    :=
    \inf
    \left\{
        0\leq t\leq T:
        M_t
        >
        \widehat q_{\alpha}(g)g_t
    \right\}.
\]

The following theorem establishes the finite-sample validity of this
calibrated boundary. The argument is purely rank-based and therefore does
not require knowledge of the null distribution of the e-process.

\begin{theorem}[Finite-sample validity of reference-null calibrated
e-process boundaries]
\label{thm:eprocess-boundary}

Suppose that, under \(P_0\), the test trajectory and the \(B\) independent
calibration trajectories are exchangeable.
Then the stopping rule defined above satisfies
\(P_0(\nu_\alpha\leq T) \leq \alpha .\)

\end{theorem}

\begin{remark}
The theorem shows that reference-null calibration provides a
distribution-specific replacement for the Ville boundary. The validity
does not follow from a sharper version of Ville's inequality; instead, it
is obtained by calibrating the finite-horizon null distribution of the
path statistic directly. In particular, the result remains valid even
when the null crossing distribution is analytically intractable.
\end{remark}

The role of the reference trajectories is therefore analogous to their
role in conformal inference: they provide an exchangeable reference
distribution for a statistic whose exact null law is unknown. The price
for this improvement is that the resulting boundary is horizon-specific
and requires independent null trajectories, whereas the universal
boundary \(1/\alpha\) applies without additional information and over an
arbitrary monitoring horizon.
Algorithm~\ref{alg:reference-null-eprocess} summarizes the complete calibration and monitoring procedure.

\begin{algorithm}[t]
\caption{Reference-null calibration of a finite-horizon e-process}
\label{alg:reference-null-eprocess}
\begin{algorithmic}[1]

\Require
Significance level \(\alpha\in(0,1)\), monitoring horizon \(T\),
number of reference-null trajectories \(B\), deterministic positive
boundary template \(g=(g_0,\ldots,g_T)\), and an e-process construction
\(\mathcal M\).

\State Generate \(B\) independent reference-null data trajectories
under \(P_0\).

\For{$b=1,\ldots,B$}

    \State Apply the same e-process construction \(\mathcal M\) to the
    \(b\)-th reference trajectory to obtain
    \(M_0^{\mathrm{ref},b}, M_1^{\mathrm{ref},b}, \ldots, M_T^{\mathrm{ref},b}.\)

    \State Compute its normalized path maximum
    \(R_g^{(b)} \gets \max_{0\leq t\leq T} \frac{M_t^{\mathrm{ref},b}}{g_t}.\)

\EndFor

\State Sort the reference statistics:
\(R_{g,(1)} \leq\cdots\leq R_{g,(B)}.\)

\State Set
\(k_{\alpha,B} \gets \left\lceil (1-\alpha)(B+1) \right\rceil.\)

\If{$k_{\alpha,B}=B+1$}

    \State Set
    \(\widehat q_\alpha(g)\gets+\infty.\)

\Else

    \State Set
    \(\widehat q_\alpha(g) \gets R_{g,(k_{\alpha,B})}.\)

\EndIf

\State Define the calibrated boundary
\(\widehat b_t \gets \widehat q_\alpha(g)g_t, \qquad 0\leq t\leq T.\)

\State Initialize the test e-process with
\(M_0^{\mathrm{test}}\gets1.\)

\For{$t=0,1,\ldots,T$}

    \If{$t\geq1$}
        \State Observe the new test data and update
        \(M_t^{\mathrm{test}}\) according to \(\mathcal M\).
    \EndIf

    \If{$M_t^{\mathrm{test}}>\widehat b_t$}

        \State \Return
        \(\nu_\alpha\gets t.\)

    \EndIf

\EndFor

\State \Return
\(\nu_\alpha\gets\infty.\)

\end{algorithmic}
\end{algorithm}

To quantify the possible reduction relative to the Ville boundary, consider
the oracle upper quantile of the null path maximum. Define
\[
    s_{\gamma,T}
    :=
    \inf
    \left\{
        s>0:
        P_0(S_T\geq s)\leq\gamma
    \right\},
    \qquad
    S_T=\max_{0\leq t\leq T}M_t .
\]
The quantity \(s_{\gamma,T}\) represents the smallest
finite-horizon level-\(\gamma\) boundary for the particular e-process,
whereas \(1/\alpha\) is the universal level-\(\alpha\) boundary obtained
from Ville's inequality. The following proposition gives a finite-sample
upper bound on the calibrated boundary in terms of this oracle quantity.

\begin{theorem}[High-probability reduction of the calibrated boundary]
\label{thm:boundary-reduction}

Fix a finite horizon \(T\) and a deterministic positive boundary template
\(g=(g_0,\ldots,g_T)\). 
Define the aspect ratio
\(K_g := \frac{ \max_{0\leq t\leq T}g_t }{ \min_{0\leq t\leq T}g_t } .\)
For \(0<\gamma<1\) and \(\eta>0\),
\[
P
\left\{
    \max_{0\leq t\leq T}
    \widehat q_\alpha(g)g_t
    <
    K_gs_{\alpha-\eta,T}
\right\}
\geq
1-
\exp
\left[
-2B
\left(
\eta-\frac{1-\alpha}{B}
\right)_+^2
\right].
\]
Consequently, with the same probability,
\[
\frac{
    \max_{0\leq t\leq T}
    \widehat q_\alpha(g)g_t
}{
    1/\alpha
}
\leq
\alpha K_gs_{\alpha-\eta,T}.
\]

\end{theorem}

\begin{remark}
The reduction in Theorem~\ref{thm:boundary-reduction} is governed by
two distinct factors. First, consider the constant boundary template
\(g_t\equiv1\), for which \(K_g=1\). In this case, the proposition gives
$
    \frac{
        \widehat q_\alpha(g)
    }{
        1/\alpha
    }
    \leq
    \alpha s_{\alpha-\eta,T}
$
with high probability. Hence, whenever the actual finite-horizon null
quantile of the path maximum is strictly smaller than the Ville bound,
the reference-null calibrated boundary improves upon \(1/\alpha\). In
particular, as the number of calibration trajectories \(B\) increases, the
term
$
    \frac{1-\alpha}{B}
$
becomes negligible and \(\eta\) can be chosen arbitrarily small. Therefore,
the calibrated boundary approaches the oracle finite-horizon quantile and
can be strictly smaller than the universal Ville boundary whenever the
null distribution of the e-process is not least favorable.

The second factor is introduced by a nonconstant boundary template. The
aspect ratio
\(K_g = \frac{\max_{0\leq t\leq T}g_t} {\min_{0\leq t\leq T}g_t}\)
quantifies the price paid for allowing temporal variation in the boundary.
To see this point, consider a simple two-step example with
$
    g_1=\varepsilon,g_2=1 ,
$
where \(\varepsilon>0\) is small. Then
$
    K_g=\frac{1}{\varepsilon}.
$
Suppose that under the null hypothesis the e-process is the constant
process
$
    M_1=M_2=1 .
$
For this process,
$
    R_g
    =
    \max\left\{
        \frac{M_1}{g_1},
        \frac{M_2}{g_2}
    \right\}
    =
    \frac1{\varepsilon}.
$
Hence the calibrated multiplier must be of order
\(1/\varepsilon\), and the resulting boundary at the second time point is
also inflated to the same order$
    \widehat q_\alpha(g)g_2
    \approx
    \frac1{\varepsilon}.
$
Therefore, an unnecessarily small value of \(g_t\) at one time point can
force the entire calibrated boundary to increase, even though the
underlying e-process itself contains no evidence against the null.

Consequently, without additional prior knowledge about the temporal
behavior of the e-process, the natural choice is the constant template
\(g_t\equiv1\). All subsequent developments in this paper adopt this
choice and calibrate the path maximum directly.
\end{remark}

The previous propositions establish validity of the reference-null
calibrated boundary and characterize when it can improve upon the
universal Ville boundary. We next quantify the finite-sample accuracy of
the empirical calibration procedure. In particular, we compare the
empirical multiplier \(\widehat q_\alpha(g)\) with the oracle multiplier
that would be obtained if the null distribution of the normalized path
maximum were known.

\begin{proposition}[Accuracy of the empirical reference-null boundary]
\label{prop:boundary-accuracy}

Let
\[
    q_{\alpha,T}^{\star}(g)
    :=
    \inf
    \left\{
        q>0:
        \mathbb P_0(R_g\leq q)
        \geq
        1-\alpha
    \right\}
\]
denote the oracle finite-horizon calibration multiplier. For any
\(\delta\in(0,1)\), define
\(\varepsilon_B(\delta) := \sqrt{ \frac{\log(2/\delta)}{2B} } + \frac{2}{B}.\)
Then, with probability at least \(1-\delta\),
\[
    q_{\alpha+\varepsilon_B(\delta),T}^{\star}(g)
    \leq
    \widehat q_{\alpha}(g)
    \leq
    q_{\alpha-\varepsilon_B(\delta),T}^{\star}(g),
\]
provided that
$
    \alpha>\varepsilon_B(\delta).
$

\end{proposition}

\begin{remark}
Proposition~\ref{prop:boundary-accuracy} is a direct consequence of a
DKW-type uniform deviation bound for the empirical distribution of
\(R_g\). It shows that the empirical reference-null boundary consistently
approximates the oracle finite-horizon boundary, with the quantile error
controlled at the usual rate \(O(B^{-1/2})\).
\end{remark}

\section{Conformal Martingales for Reference-Null Calibration}
\label{sec:conformal-martingale}

We now specialize the reference-null calibration framework to adaptive
conformal martingales. These processes provide a useful setting for
studying calibration gains: their null law can be simulated exactly,
while their adaptive betting rules allow us to analyze evidence
accumulation under distributional changes. Building on conformal
martingales and histogram betting
\citep{vovk2003testing,vovk2005algorithmic,
fedorova2012plugin,eliades2020histogram}, we specify a
Krichevsky--Trofimov (KT) histogram process and its restart mixture.
These constructions provide the basis for the finite-sample power and
detection-delay analysis in Section~\ref{sec:power-delay}, where we
quantify how calibrated boundaries interact with learning and
restart aggregation.

\subsection{Conformal martingales and their pivotal null law}
\label{subsec:conformal-construction}

Let \(X_1,X_2,\ldots\) be real-valued scores and consider the
exchangeability null
\[
    \mathcal P_0
    :=
    \left\{
        P:(X_t)_{t\geq1}\text{ is exchangeable under }P
    \right\}.
\]
The scores may be obtained from more complex observations, provided that
the resulting score sequence is exchangeable under the null.
Let \(U_1,U_2,\ldots\) be independent \(\mathrm{Unif}(0,1)\)
random variables, independent of the entire score sequence. Define
the smoothed conformal \(p\)-values by
\begin{equation}
\label{eq:section3-conformal-p}
    p_t
    :=
    \frac{
        \sum_{i=1}^{t}\mathbf 1\{X_i>X_t\}
        +
        U_t\sum_{i=1}^{t}\mathbf 1\{X_i=X_t\}
    }{t},
    \qquad t\geq1.
\end{equation}
A standard randomized conformal validity result implies that, under
every \(P\in\mathcal P_0\), the sequential conformal \(p\)-values are
independent and uniformly distributed on \([0,1]\);
see \citet{vovk2003testing,vovk2005algorithmic,vovk2019nonparametric}.
Write
$
    \mathcal G_t:=\sigma(p_1,\ldots,p_t),
    \mathcal G_0:=\{\varnothing,\Omega\}.
$
For each \(t\geq1\), let \(f_t:[0,1]\to[0,\infty)\) be a
\(\mathcal G_{t-1}\)-measurable random function satisfying
$
\int_0^1 f_t(u)\,\mathrm du=1
$
almost surely.  The associated
conformal martingale is
\[
    M_0:=1,
    \qquad
    M_t:=\prod_{i=1}^{t}f_i(p_i).
\]
Indeed, under every \(P\in\mathcal P_0\),
$
    \mathbb E_P[M_t\mid\mathcal G_{t-1}]
    =
    M_{t-1}\int_0^1f_t(u)\,\mathrm du
    =
    M_{t-1}.
$
Thus \(M_t\) is a nonnegative martingale, and hence an e-process, for
the exchangeability null with respect to the \(p\)-value filtration.
In particular,
\[
    \sup_{P\in\mathcal P_0}
    P\left(
        \sup_{t\geq0}M_t\geq\frac1\alpha
    \right)
    \leq\alpha.
\]
All betting rules considered below use the observations only through
the past conformal \(p\)-values. Since the null distribution of the
conformal \(p\)-value sequence is invariant over
\(P\in\mathcal P_0\), the resulting martingale trajectory distribution
is also invariant over the null class. In particular, once the betting
rules are specified, their null trajectory distribution is identical for
all \(P\in\mathcal P_0\) and can be generated by replacing the conformal
\(p\)-values with independent \(\mathrm{Unif}(0,1)\) random variables.

\subsection{Adaptive KT histogram betting}
\label{subsec:KT-histogram}

Under an alternative distribution, the conformal \(p\)-values need not remain
independent or identically distributed. Let
\(q_t(\cdot\mid\mathcal G_{t-1})\) denote the conditional density of \(p_t\)
given the preceding conformal \(p\)-values. For any predictable betting density
\(f_t\), the conditional expected log growth satisfies
\[
\mathbb E[\log f_t(p_t)\mid\mathcal G_{t-1}]
=
\int_0^1 q_t(u\mid\mathcal G_{t-1})
\log q_t(u\mid\mathcal G_{t-1})\,\mathrm du
-
D_{\mathrm{KL}}
\left(
q_t(\cdot\mid\mathcal G_{t-1})
\|f_t
\right).
\]
Hence the conditionally log-optimal betting rule is
\(f_t^\star(\cdot)=q_t(\cdot\mid\mathcal G_{t-1})\). This conditional
formulation is important because, after a distributional change, the law of
\(p_t\) generally depends on the preceding conformal \(p\)-values.

Histogram-based betting provides a particularly simple adaptive approximation
to this oracle rule and has previously been used as
computationally convenient betting functions for conformal martingales
\citep{eliades2020histogram}. In the present setting, there is also a
structural reason for considering piecewise-constant betting functions: the
conditional density of a randomized conformal \(p\)-value, given its preceding
conformal \(p\)-values, is itself piecewise constant on the regular grid with
mesh size \(1/t\).

To make this observation precise, define
\(R_t^>:=\sum_{i=1}^t\mathbf 1\{X_i>X_t\}\) and
\(R_t^=:=\sum_{i=1}^t\mathbf 1\{X_i=X_t\}\). Then the randomized conformal
\(p\)-value in \eqref{eq:section3-conformal-p} can be written as
\(p_t=(R_t^>+U_tR_t^=)/t\).

\begin{proposition}[Piecewise-constant conditional density of a randomized conformal \(p\)-value]
\label{prop:piecewise-conformal-density}
For each fixed \(t\geq1\), suppose that \(U_t\sim\Unif(0,1)\) is independent
of \(X_1,\ldots,X_t\) and of the previous randomizers. Then, conditional on
\(\mathcal G_{t-1}\), the randomized conformal \(p\)-value \(p_t\) admits a
conditional density \(q_t(u\mid\mathcal G_{t-1})\) on \((0,1)\) satisfying
\begin{equation}
\label{eq:piecewise-conformal-density}
q_t(u\mid\mathcal G_{t-1})
=
\mathbb E
\left[
\left.
\frac{t}{R_t^=}
\mathbf 1
\left\{
\frac{R_t^>}{t}
<
u
<
\frac{R_t^>+R_t^=}{t}
\right\}
\,\right|\,
\mathcal G_{t-1}
\right].
\end{equation}
Consequently, \(q_t(\cdot\mid\mathcal G_{t-1})\) is constant on each interval
\(((k-1)/t,k/t)\), \(k=1,\ldots,t\).
\end{proposition}
Proposition~\ref{prop:piecewise-conformal-density} shows that
piecewise-constant betting functions are naturally aligned with the
conditional structure of randomized conformal \(p\)-values. Directly estimating
the full \(t\)-cell conditional density, however, would require a betting rule
whose dimension increases with time. We therefore use a fixed, coarser
partition of \([0,1]\), which provides a low-dimensional adaptive
approximation to the oracle betting density.

Fix \(J\geq2\) and partition \([0,1]\) into equal-width bins
\(I_j=[(j-1)/J,j/J)\), \(j=1,\ldots,J-1\), and
\(I_J=[(J-1)/J,1]\). Writing
\(q_{t,j}:=\mathbb P(p_t\in I_j\mid\mathcal G_{t-1})\), the conditionally
optimal betting density within this \(J\)-bin histogram class has height
\(Jq_{t,j}\) on \(I_j\).

Since the probabilities \(q_{t,j}\) are unknown, we estimate them
sequentially. Let
\(N_j(t):=\sum_{i=1}^t\mathbf 1\{p_i\in I_j\}\), with \(N_j(0)=0\), and use
the Krichevsky--Trofimov (KT) predictive probabilities
\[\widehat q_{j,t}^{\mathrm{KT}}
=\frac{N_j(t-1)+1/2}{t-1+J/2}.\]
The corresponding betting density is
\(\widehat f_t^{(J)}(u):=J\widehat q_{j,t}^{\mathrm{KT}}\) for \(u\in I_j\),
and the resulting wealth process is
\(\widehat M_0^{(J)}:=1\) and
\(\widehat M_t^{(J)}:=\prod_{i=1}^t\widehat f_i^{(J)}(p_i)\).
The half-count smoothing keeps every betting factor strictly positive,
including when a previously unobserved bin is encountered.

\begin{proposition}[Validity of KT histogram betting]
\label{prop:KT-histogram-validity}
For every \(P\in\mathcal P_0\),
\((\widehat M_t^{(J)})_{t\geq0}\) is a nonnegative
\((\mathcal G_t)\)-martingale initialized at one.
\end{proposition}

Algorithm~\ref{alg:online-unrestarted-KT} summarizes the online
implementation. Given the conformal \(p\)-value stream, the procedure only
needs to maintain the \(J\) bin counts and the current wealth.

\begin{algorithm}[t]
\caption{Online unrestarted KT histogram process}
\label{alg:online-unrestarted-KT}
\begin{algorithmic}[1]

\Require Monitoring horizon \(T\), number of bins \(J\geq2\), conformal
\(p\)-value stream \(p_1,p_2,\ldots\), and rejection boundary \(c_\alpha\).

\State Partition \([0,1]\) into equal-width bins \(I_1,\ldots,I_J\).

\State Initialize \(N_j(0)\gets0\), \(j=1,\ldots,J\), and
\(\widehat M_0^{(J)}\gets1\).

\For{$t=1,2,\ldots,T$}

    \State Observe \(p_t\) and determine \(j_t\) such that \(p_t\in I_{j_t}\).

    \State Compute the predictable KT betting factor
    \[
    \widehat f_t^{(J)}(p_t)
    \gets
    J\frac{N_{j_t}(t-1)+1/2}{t-1+J/2}.
    \]

    \State Update
    \(\widehat M_t^{(J)}
    \gets
    \widehat M_{t-1}^{(J)}
    \widehat f_t^{(J)}(p_t)\).

    \State Set \(N_{j_t}(t)\gets N_{j_t}(t-1)+1\) and
    \(N_j(t)\gets N_j(t-1)\) for \(j\neq j_t\).

    \If{\(\widehat M_t^{(J)}>c_\alpha\)}
        \State \Return detection time \(t\).
    \EndIf

\EndFor

\State \Return \(\infty\) if no crossing occurs by time \(T\).

\end{algorithmic}
\end{algorithm}

\subsection{Restart aggregation and its adaptation cost}
\label{subsec:restart-KT}

The unrestarted process above estimates its betting distribution from the
entire \(p\)-value history. Following a late change, the many
pre-change observations can delay its adaptation. To see this, suppose
that the scores are independent, with distribution \(F_0\) up to time
\(\tau\) and \(F_1\) thereafter. For \(t\geq\tau+2\), define
$
    \overline q_{j,t-1}
    :=
    \frac1{t-1-\tau}
    \sum_{r=\tau+1}^{t-1}\mathbb P(p_r\in I_j).
$
Since the pre-change \(p\)-values are uniform,
$
    \mathbb E[\widehat q_{j,t}^{\mathrm{KT}}]-\frac1J
    =
    \frac{t-1-\tau}{t-1+J/2}
    \left(\overline q_{j,t-1}-\frac1J\right).
$
The factor multiplying the average post-change departure is small
soon after a late change. This identity describes dilution in the
expected histogram probabilities and motivates restarting the
betting rule; its implications for power are established in
Section~\ref{sec:power-delay}.

Restarting and aggregating evidence processes is an established
approach to unknown change points. In particular,
\citet{saha2026nonpartitioned} aggregate point-null e-processes over
candidate starting times and take an infimum over candidate no-change
distributions. Their weighting schemes provide either average-run-length
or probability-of-false-alarm control. We use the same general
restart-and-aggregate principle for histogram betting on a conformal
\(p\)-value stream.

For each candidate restart time \(s\geq0\), set
$
    N_j^{(s)}(t)
    :=
    \sum_{i=s+1}^{t}\mathbf 1\{p_i\in I_j\}, t\geq s,
$
with an empty sum equal to zero, and define, for \(t>s\),
\[
    \widehat q_{j,t}^{(s)}
    :=
    \frac{N_j^{(s)}(t-1)+1/2}{t-s-1+J/2},
    \qquad
    \widehat f_t^{(s,J)}(u)
    :=
    J\widehat q_{j,t}^{(s)},\quad u\in I_j.
\]
The restarted component is
\begin{equation}
\label{eq:restarted-KT-component}
    \widehat M_t^{(s,J)}
    :=
    \begin{cases}
        1, & t\leq s,\\[2mm]
        \displaystyle\prod_{i=s+1}^{t}
        \widehat f_i^{(s,J)}(p_i), & t>s.
    \end{cases}
\end{equation}
Thus \(s=0\) recovers the unrestarted process. Each component resets
its betting counts and wealth but continues to use the original
\(p\)-values in \eqref{eq:section3-conformal-p}. In particular,
restart does not remove the changing composition of the reference
sample used for conformal ranking; the post-change \(p\)-values can
remain dependent and nonstationary.

Let \((\pi_s)_{s\geq0}\) be deterministic nonnegative weights with
\(\sum_{s=0}^{\infty}\pi_s=1\). Define
\begin{equation}
\label{eq:restart-mixture-process}
    \widehat M_t^{R,J}
    :=
    \sum_{s=0}^{\infty}\pi_s\widehat M_t^{(s,J)}
    =
    \sum_{s=0}^{t-1}\pi_s\widehat M_t^{(s,J)}
    +
    \sum_{s=t}^{\infty}\pi_s.
\end{equation}
The final term represents the capital allocated to components that
have not yet started. It is part of the process definition and keeps
the total initial wealth equal to one.

\begin{proposition}[Validity of restart aggregation]
\label{prop:restart-mixture-validity}
For every \(P\in\mathcal P_0\),
\((\widehat M_t^{R,J})_{t\geq0}\) is a nonnegative
\((\mathcal G_t)\)-martingale initialized at one, and hence an
e-process for \(\mathcal P_0\).
\end{proposition}

For any candidate change point \(\tau\) with \(\pi_\tau>0\),
nonnegativity gives the pathwise comparison
\begin{equation*}
\label{eq:restart-aggregation-penalty}
    \log\widehat M_t^{R,J}
    \geq
    \log\widehat M_t^{(\tau,J)}
    -
    \log\frac1{\pi_\tau}.
\end{equation*}
Thus \(\log(1/\pi_\tau)\) is an upper bound on the log-wealth
loss relative to the component started at \(\tau\). This deterministic
aggregation cost enters the power and delay analysis separately from
the KT learning cost and the rejection boundary.

\begin{remark}[Finite-horizon restart weights]
\label{thm:uniform-restart-minimax}
For a fixed monitoring horizon \(T\), one may put all weight on
\(\mathcal S_T=\{0,\ldots,T-1\}\), obtaining $
    \widehat M_t^{\mathrm{mix},J}
    :=
    \sum_{s=0}^{T-1}\pi_s\widehat M_t^{(s,J)},
    0\leq t\leq T.
$
This is the special case of \eqref{eq:restart-mixture-process}
with \(\pi_s=0\) for \(s\geq T\). Uniform weights
\begin{equation}
\label{eq:uniform-restart-weights}
    \pi_s=\frac1T,\qquad s=0,\ldots,T-1,
\end{equation}
give the same aggregation penalty \(\log T\) at every candidate
change point. They uniquely minimize
\(\max_{s\in\mathcal S_T}\log(1/\pi_s)\), since
\(\min_s\pi_s\leq1/T\), with equality only for uniform weights.
 We therefore recommend adopting \eqref{eq:uniform-restart-weights} in the absence of prior
information
\end{remark}

Algorithm~\ref{alg} summarizes the online
implementation. At time \(t\), only the \(t\) started components need updating;
the remaining components are represented by the weight tail in
\eqref{eq:restart-mixture-process}. Storing previous bin labels and
scanning them in reverse order gives the counts needed to update all
active components in \(O(t)\) time per observation and \(O(t)\)
memory, in addition to the cost of producing the conformal
\(p\)-values. Thus the direct implementation takes \(O(T^2)\) time
over a horizon \(T\).
\begin{algorithm}[t]
\caption{Online restart-mixture KT histogram process}
\label{alg}
\begin{algorithmic}[1]

\Require
Number of bins \(J\geq1\), restart weights
\(\{\pi_s\}_{s\geq0}\) satisfying
\(\pi_s\geq0,\qquad \sum_{s=0}^{\infty}\pi_s=1,\)
a sequential conformal \(p\)-value stream
\(p_1,p_2,\ldots\), and rejection boundary \(c_\alpha\).

\State Partition \([0,1]\) into equal-width bins
\(I_1,\ldots,I_J\).

\State Initialize
\(\overline\Pi_0 \gets \sum_{s=0}^{\infty}\pi_s=1 .\)

\State Initialize empty arrays for the historical bin labels
\(\{j_\ell\}\) and restart component wealths
\(\{\widehat M^{(s,J)}\}\).

\For{$t=1,2,\ldots$}

\State Observe \(p_t\) and determine \(j_t\).

\State Start the new restart component
\(\widehat M_{t-1}^{(t-1,J)}\gets1 .\)

\State Update the restart-weight tail:
\(\overline\Pi_t \gets \overline\Pi_{t-1}-\pi_{t-1}.\)

\State Set
\(C\gets0, \qquad S_t\gets0 .\)

\For{$s=t-1,t-2,\ldots,0$}

\State Update the \(s\)-th component:
\(\widehat M_t^{(s,J)} \gets \widehat M_{t-1}^{(s,J)} J \frac{ C+1/2 } { t-s-1+J/2 }.\)

\State Accumulate:
\(S_t \gets S_t+ \pi_s \widehat M_t^{(s,J)} .\)

\If{$s\geq1$}

\State Update
\(C \gets C+\mathbf 1\{j_s=j_t\}.\)

\EndIf

\EndFor

\State Form the restart mixture:
\(\widehat M_t^{R,J} \gets S_t+\overline\Pi_t .\)

\State Store current bin label \(j_t\).

\If{$\widehat M_t^{R,J}\geq c_\alpha$}

\State \Return detection time \(t\).

\EndIf

\EndFor

\end{algorithmic}
\end{algorithm}

\subsection{Applying reference-null calibration}
\label{subsec:conformal-reference-calibration}

Fix \(T\), \(J\), the restart weights when applicable, and either of
the preceding process constructions, denoted by \(M\).
Generate \(B\) independent reference sequences of \(T\) i.i.d.\
\(\mathrm{Unif}(0,1)\) variables, independently of the test data,
and apply the same construction to each sequence. With
\(S_T^{(b)}:=\max_{0\leq t\leq T}M_t^{(b)}\), set
\[
    \widehat c_{\alpha,B,T}
    :=
    S_{T,(k_{\alpha,B})},
    \qquad
    k_{\alpha,B}:=\left\lceil(1-\alpha)(B+1)\right\rceil,
\]
where \(S_{T,(1)}\leq\cdots\leq S_{T,(B)}\) are the ordered
reference maxima and \(S_{T,(B+1)}:=+\infty\). All design choices
are fixed independently of the calibration sample.

The pivotal null law and Theorem~\ref{thm:eprocess-boundary} give
\[
    \sup_{P\in\mathcal P_0}
    P\left(
        \max_{0\leq t\leq T}M_t
        >
        \widehat c_{\alpha,B,T}
    \right)
    \leq\alpha,
\]
where probability includes calibration randomness. This is the
constant-template application of Section~\ref{subsec:eprocess-calibration}.
Each process is calibrated using its own reference maxima.
Calibration leaves its betting and restart rules unchanged.
Section~\ref{sec:power-delay} quantifies the power and delay gains
when its calibrated boundary is smaller than \(1/\alpha\).

\section{Power and Detection-Delay Analysis}
\label{sec:power-delay}

The previous sections separate two ingredients of sequential detection.
The betting strategy determines how rapidly the e-process accumulates
evidence under a distributional change, whereas the rejection boundary
determines how much evidence is required before detection.  We now study
these two ingredients jointly.  Our main goals are to characterize the
power of the conformal likelihood-ratio benchmark and the practical KT
histogram procedures, to quantify the effect of restart aggregation, and
to show explicitly how a reduction of the rejection boundary translates
into gains in power and detection delay.

The reference-null calibrated boundary is random because it is computed
from an independent calibration sample.  Throughout this section, we first
condition on this calibration sample.  Conditional on the calibration data,
the resulting boundary is fixed and independent of the test trajectory, so
the power and delay results below may be applied exactly as for a
deterministic boundary.  Unconditional statements can then be recovered by
averaging over the calibration sample.

We work under the change-point model
\begin{equation}
\label{eq:section4-change-point-model}
    X_1,\ldots,X_\tau
    \overset{\mathrm{iid}}{\sim}F_0,
    \qquad
    X_{\tau+1},X_{\tau+2},\ldots
    \overset{\mathrm{iid}}{\sim}F_1,
\end{equation}
where the pre-change and post-change samples are independent and
\(F_0,F_1\) are continuous distributions.  The corresponding smoothed
conformal \(p\)-values are
$$
    p_t
    =
    \frac{
        \displaystyle
        \sum_{r=1}^{t-1}
        \mathbf 1\{X_r>X_t\}
        +
        U_t
    }{t},
    \qquad
    U_t\overset{\mathrm{iid}}{\sim}\operatorname{Unif}(0,1),
$$
where the randomizers are independent of the scores.

A convenient measure of the separation between the pre-change and
post-change score distributions is
\begin{equation}
\label{eq:AUC-separation}
    \theta
    :=
    \mathbb P(X_0>X_1),
    \qquad
    \Delta
    :=
    \left|
        \theta-\frac12
    \right|,
\end{equation}
where \(X_0\sim F_0\) and \(X_1\sim F_1\) are independent.  The quantity
\(\theta\) is the usual AUC-type pairwise ranking probability.  Under the
absence of a distributional change, \(F_0=F_1\) and
\(\theta=1/2\).  Hence \(\Delta>0\) measures a directional departure from
the exchangeability null directly at the level of the underlying score distributions.

\subsection{An oracle conformal likelihood-ratio benchmark}
\label{subsec:oracle-power}

We first consider an oracle process that uses the exact conditional law of
the conformal \(p\)-value sequence under the alternative.  This provides a
benchmark for the amount of evidence that can in principle be extracted
from the conformal \(p\)-values when the post-change distribution is
known.

Fix a change point \(\tau\geq1\) and a post-change monitoring length
\(n\geq1\), and let
\(T:=\tau+n.\)
For \(1\leq t\leq T\), let \(P_t\) denote the joint law of
\(p_{1:t}\) under the change-point model
\eqref{eq:section4-change-point-model}, and let
\(Q_t := \operatorname{Unif}(0,1)^{\otimes t}.\)
Let
\(q_t \left( u\mid p_{1:t-1} \right)\)
denote the conditional density of \(p_t\) under the alternative.
As discussed in Section~\ref{subsec:KT-histogram}, this conditional
density is the log-optimal predictable betting function when the
alternative law is known.

Define the oracle conformal likelihood-ratio process by
\begin{equation}
\label{eq:oracle-conformal-LR}
    M_t^\star
    :=
    \prod_{s=1}^t
    q_s
    \left(
        p_s\mid p_{1:s-1}
    \right),
    \qquad
    M_0^\star=1.
\end{equation}
By the chain rule for conditional densities,
\(M_t^\star = \frac{dP_t}{dQ_t}, \qquad 1\leq t\leq T.\)
Consequently, \((M_t^\star)_{t\leq T}\) is a nonnegative martingale
under \(Q_T\).

For a deterministic positive boundary
\((B_t)_{1\leq t\leq T}\), define
\(D^\star := \inf \left\{ 1\leq t\leq T: M_t^\star>B_t \right\},\)
with \(\inf\varnothing=\infty\).

\begin{theorem}[Finite-sample power bound for the oracle conformal
likelihood-ratio process]
\label{thm:oracle-conformal-power}

Suppose that the AUC separation defined in
\eqref{eq:AUC-separation} satisfies \(\Delta>0\), and define
\(r_{\tau,n} := \frac{\tau n}{\tau+n}.\)
For \(0<\lambda<1\), let
$$
    c_\lambda
    :=
    \frac{
        2\lambda(1-\lambda)
    }{
        \left(
            \sqrt{\lambda}
            +
            \sqrt{1-\lambda}
        \right)^2
    }.
$$
Then
\begin{equation}
    \mathbb P_\tau
    \left(
        D^\star>T
    \right)
    \leq
    \inf_{0<\lambda<1}
    2
    \exp
    \left\{
        \lambda\log B_T
        -
        c_\lambda
        \Delta^2 r_{\tau,n}
    \right\}.
\label{eq:oracle-power-bound-general}
\end{equation}

In particular, taking \(\lambda=1/2\) gives
\begin{equation}
    \mathbb P_\tau
    \left(
        D^\star>T
    \right)
    \leq
    2B_T^{1/2}
    \exp
    \left\{
        -\frac{\Delta^2}{4}
        \frac{\tau n}{\tau+n}
    \right\}.
\label{eq:oracle-power-bound}
\end{equation}
\end{theorem}

\begin{remark}
\label{rem:oracle-power-interpretation}

Theorem~\ref{thm:oracle-conformal-power} separates the difficulty of the
change-point problem into a signal term and a boundary term.  The signal
is governed by
$
    \Delta^2 r_{\tau,n}
$
where
\(\frac12\min\{\tau,n\} \leq r_{\tau,n} \leq \min\{\tau,n\}.\)
Thus \(r_{\tau,n}\) has the same order as the smaller of the pre-change
and post-change sample sizes.

For any fixed \(\lambda\in(0,1)\), the general bound
\eqref{eq:oracle-power-bound-general} implies that
\[
    c_\lambda
    \Delta^2
    \frac{\tau n}{\tau+n}
    -
    \lambda\log B_T
    \longrightarrow\infty
\]
is sufficient for
\(\mathbb P_\tau(D^\star\leq T) \longrightarrow1.\)
For the constant Ville boundary
\(B_t=\frac1\alpha,\)
with fixed \(\alpha\in(0,1)\), this reduces to
\(\Delta^2 \frac{\tau n}{\tau+n} \longrightarrow\infty.\)
Hence, if \(F_0\) and \(F_1\) are fixed and \(\Delta>0\), then
\(\min\{\tau,n\}\longrightarrow\infty\)
is sufficient for the oracle detection probability to converge to one.

The general form
\eqref{eq:oracle-power-bound-general} also gives a direct quantitative
description of the effect of reducing the rejection boundary.  Suppose,
for example, that
\(B_T^{\mathrm C} = \rho B_T^{\mathrm V}, \qquad 0<\rho<1,\)
where \(B_T^{\mathrm V}=1/\alpha\) is the Ville boundary.  For every fixed
\(\lambda\in(0,1)\), the corresponding upper bound on the non-detection
probability is reduced by the multiplicative factor
$
    \rho^\lambda.
$

Equivalently, for a prescribed target power \(1-\beta\), a sufficient
condition obtained from any fixed \(\lambda\in(0,1)\) is
\begin{equation}
\label{eq:oracle-target-power-condition}
    c_\lambda
    \Delta^2
    \frac{\tau n}{\tau+n}
    \geq
    \lambda\log B_T
    +
    \log\frac{2}{\beta}.
\end{equation}
Thus replacing \(B_T^{\mathrm V}\) by
\(B_T^{\mathrm C}=\rho B_T^{\mathrm V}\) reduces the required right-hand
side by
$
    \lambda\log\frac1\rho.
$
This makes explicit how a smaller calibrated boundary translates into a
weaker signal requirement for the same guaranteed power.

The proof of Theorem~\ref{thm:oracle-conformal-power} is based on a
rank-separation argument.  Although the oracle process is defined through
the joint likelihood ratio of the conformal \(p\)-value sequence, the
randomized conformal \(p\)-values retain the sequential insertion ranks of
the original observations.  Indeed, under the continuity assumption,
\[
    p_t
    =
    \frac{J_t+U_t}{t},
    \qquad
    J_t
    :=
    \sum_{s=1}^{t-1}
    \mathbf 1\{X_s>X_t\},
\]
with \(0<U_t<1\) almost surely, and hence
\(J_t=\lfloor tp_t\rfloor .\)
The sequence \((J_1,\ldots,J_T)\) is an insertion-rank representation of
the complete relative ordering of
\((X_1,\ldots,X_T)\): recursively, once the ordering of the first
\(t-1\) observations is known, \(J_t\) determines the unique position at
which \(X_t\) is inserted, since exactly \(J_t\) preceding observations
are larger than \(X_t\).  Consequently, every pairwise comparison
\(\mathbf 1\{X_i>X_j\}\), and therefore the two-sample Mann--Whitney
statistic
\[
    \widehat\theta
    :=
    \frac1{\tau n}
    \sum_{i=1}^{\tau}
    \sum_{j=1}^{n}
    \mathbf 1\{X_i>X_{\tau+j}\},
\]
is measurable with respect to the conformal \(p\)-value sequence.

This observation allows the separation between the change-point law and
the product-uniform null law to be studied directly through
\(\widehat\theta\).  Under the alternative,
\(\mathbb E_\tau[\widehat\theta] = \theta = \mathbb P(X_0>X_1),\)
whereas under the product-uniform null law the induced ranks are those of
i.i.d.\ continuous observations and
\(\mathbb E_{Q_T}[\widehat\theta] = \frac12.\)
Thus the same rank-measurable statistic separates the two joint
\(p\)-value laws by
\(\Delta = \left| \theta-\frac12 \right|.\)
A concentration inequality for \(\widehat\theta\) then yields an event
having high probability under the change-point law and low probability
under the null law.  This event controls the Hellinger, or more generally
R\'enyi, overlap between the two joint \(p\)-value distributions, which is
subsequently converted into a finite-sample bound on the probability that
the oracle likelihood-ratio process fails to cross the rejection
boundary.

This proof structure also explains the appearance of
\(\Delta^2 \frac{\tau n}{\tau+n}\)
in the resulting power bound.  The statistic \(\widehat\theta\) is the
classical two-sample Mann--Whitney U-statistic estimating
\(\theta = \mathbb P(X_0>X_1).\)
Under the usual nondegeneracy conditions, its Hoeffding decomposition
implies
\[
    \operatorname{Var}(\widehat\theta)
    =
    \frac{\sigma_0^2}{\tau}
    +
    \frac{\sigma_1^2}{n}
    +
    O\left(\frac1{\tau n}\right),
\]
for finite constants \(\sigma_0^2\) and \(\sigma_1^2\);
see \citet{mann1947test,hoeffding1948class} and the general rank-test
theory in \citet{hajek1999theory}.  Hence, in the nondegenerate case,
the stochastic fluctuation of \(\widehat\theta\) is of order
\[
    \left(
        \frac1\tau+\frac1n
    \right)^{1/2}
    =
    \left(
        \frac{\tau n}{\tau+n}
    \right)^{-1/2}.
\]
Therefore
\(\Delta \sqrt{ \frac{\tau n}{\tau+n} }\)
is the natural standardized signal strength for rank-based separation.
Accordingly, alternatives satisfying
\(\Delta = O\left[ \left( \frac{\tau n}{\tau+n} \right)^{-1/2} \right]\)
lie on the usual local scale for two-sample rank inference, whereas
\(\Delta^2 \frac{\tau n}{\tau+n} \longrightarrow\infty\)
places the AUC separation above this stochastic scale and is therefore a
natural sufficient regime for the oracle detection probability to
converge to one.

Theorem~\ref{thm:oracle-conformal-power} provides an explicit
finite-sample guarantee in terms of the interpretable AUC separation
\(\Delta\).  Although the oracle likelihood-ratio process is defined
through the full joint law of the conformal \(p\)-value sequence, the
result shows that a single rank-based separation measure is already
sufficient to guarantee exponential decay of the non-detection
probability.  In particular, the bound depends on the alternative only
through
$
    \Delta^2\frac{\tau n}{\tau+n},
$
which yields a transparent connection between distributional separation,
effective sample size, and detection power.

\end{remark}

\subsection{Finite-sample power of KT histogram processes}
\label{subsec:KT-power}

We next study the power of the implementable KT histogram procedures
introduced in Section~\ref{sec:conformal-martingale}.  Unlike the oracle
likelihood-ratio process, the KT procedures do not know the conditional
alternative law of the conformal \(p\)-values and must learn a betting
distribution sequentially.  Their power therefore reflects three
distinct effects: the strength of the post-change rank signal, the
approximation and learning cost of histogram betting, and, for the
restart mixture, the additional cost of adapting to an unknown change
point.

The two KT procedures admit a common proof strategy.  First, the KT
predictive distribution achieves, up to logarithmic regret, the wealth
of the best fixed \(J\)-bin histogram fitted to the observed
\(p\)-values.  Second, the empirical histogram divergence from the
uniform distribution can be lower bounded by the deviation of the
average conformal \(p\)-value from \(1/2\).  Finally, this average has an
explicit mean under the change-point model and satisfies a
bounded-difference concentration inequality.  Combining these three
ingredients yields the finite-sample power bounds below.

Let
\(\theta := \mathbb P(X_0>X_1), \qquad \Delta := \left| \theta-\frac12 \right|,\)
where \(X_0\sim F_0\) and \(X_1\sim F_1\) are independent, and define
\(\varepsilon := \operatorname{sign} \left( \theta-\frac12 \right).\)
Throughout this subsection we assume \(\Delta>0\).

\subsubsection{A common pathwise lower bound for KT wealth}
\label{subsubsec:KT-pathwise-bound}

We first record the deterministic ingredient underlying both power
results.  Consider any block \(u_1,\ldots,u_N\in[0,1]\), and let
\[
    N_j
    :=
    \sum_{\ell=1}^N
    \mathbf 1\{u_\ell\in I_j\},
    \qquad
    \widehat q_j
    :=
    \frac{N_j}{N},
\]
with
\[
    \widehat{\boldsymbol q}
    :=
    (\widehat q_1,\ldots,\widehat q_J),
    \qquad
    \boldsymbol u_J
    :=
    \left(
        \frac1J,\ldots,\frac1J
    \right).
\]

For a generic block \(u_1,\ldots,u_N\), let

$$
    M_N^{\mathrm{KT}}(u_{1:N})
    :=
    \prod_{\ell=1}^N
    J
    \frac{
        N_{\ell-1,j(u_\ell)}+1/2
    }{
        \ell-1+J/2
    },
$$

where \(j(u_\ell)\) is the index of the bin containing \(u_\ell\), and

\(N_{\ell-1,j} := \sum_{r=1}^{\ell-1} \mathbf 1\{u_r\in I_j\}.\)

Thus \(M_N^{\mathrm{KT}}\) is exactly the KT betting product applied to
the block \(u_{1:N}\).  In particular,

$$
    M_m^{\mathrm{KT}}
    \bigl(
        p_{\tau+1:\tau+m}
    \bigr)
    =
    \widehat M_{\tau+m}^{(\tau,J)},
    \qquad
    M_{\tau+m}^{\mathrm{KT}}
    \bigl(
        p_{1:\tau+m}
    \bigr)
    =
    \widehat M_{\tau+m}^{(J)}.
$$

Hence the restarted and unrestarted procedures differ only in the block
of conformal \(p\)-values to which the same KT predictor is applied.

If \(N_j\) denotes the total number of observations in bin \(I_j\) over
the block \(u_{1:N}\), then the KT predictive product admits the standard
Dirichlet-\(1/2\) representation

$$
    M_N^{\mathrm{KT}}
    =
    J^N
    \frac{\Gamma(J/2)}
         {\Gamma(N+J/2)}
    \prod_{j=1}^J
    \frac{\Gamma(N_j+1/2)}
         {\Gamma(1/2)}.
$$

Its logarithm is within \(O(\log N)\) of the best fixed histogram
likelihood.  More precisely, standard Stirling bounds imply the following
pathwise result; see also the universal-coding interpretation of the KT
predictive distribution in
\citet{krichevsky1981performance}.

\begin{lemma}[Pathwise KT evidence bound]
\label{lem:KT-pathwise-evidence}

For every fixed \(J\geq2\), there exists a constant \(C_J\geq0\),
depending only on \(J\), such that for every \(N\geq1\),
$$
    \log M_N^{\mathrm{KT}}
    \geq
    N
    D_{\mathrm{KL}}
    \left(
        \widehat{\boldsymbol q}
        \,\middle\|\,
        \boldsymbol u_J
    \right)
    -
    \frac{J-1}{2}\log(N+1)
    -
    C_J.
$$
Moreover, if
\(\overline u_N := \frac1N\sum_{\ell=1}^N u_\ell,\)
then
\begin{equation}
\label{eq:histogram-KL-average-p}
    D_{\mathrm{KL}}
    \left(
        \widehat{\boldsymbol q}
        \,\middle\|\,
        \boldsymbol u_J
    \right)
    \geq
    2
    \left(
        \left|
            \overline u_N-\frac12
        \right|
        -
        \frac1{2J}
    \right)_+^2.
\end{equation}
Consequently,
$$
\log M_N^{\mathrm{KT}}
\geq
2N
\left(
    \left|
        \overline u_N-\frac12
    \right|
    -
    \frac1{2J}
\right)_+^2
-
\frac{J-1}{2}\log(N+1)
-
C_J.
$$
\end{lemma}

The term
$
    \frac1{2J}
$
has a simple origin.  Replacing each \(p\)-value by the midpoint of its
histogram bin changes it by at most \(1/(2J)\).  The empirical mean of
these midpoints is a bounded linear functional of the empirical
histogram, and Pinsker's inequality then converts its departure from
\(1/2\) into the KL lower bound
\eqref{eq:histogram-KL-average-p}.  Thus \(1/(2J)\) is a discretization
cost, while
\(\frac{J-1}{2}\log(N+1)+C_J\)
is the KT learning regret.

\subsubsection{Restart-mixture KT process}
\label{subsubsec:restart-KT-power}

We first analyze the restart mixture.  For \(m\geq1\), define the
post-change average
\[
    \overline p_{\tau,m}
    :=
    \frac1m
    \sum_{k=1}^m
    p_{\tau+k},
\]
and let
\begin{equation}
    a_{\tau,m}
    :=
    \frac1m
    \sum_{k=1}^m
    \frac{\tau}{\tau+k}.
    \label{eq:a-tau-m}
\end{equation}

A direct calculation under the change-point model gives
$$
    \varepsilon
    \left(
        \mathbb E_\tau
        \overline p_{\tau,m}
        -
        \frac12
    \right)
    =
    \Delta a_{\tau,m}.
$$
Indeed, for every \(k\geq1\),
\[
    \mathbb E_\tau[p_{\tau+k}]
    =
    \frac12
    +
    \frac{\tau}{\tau+k}
    \left(
        \theta-\frac12
    \right).
\]
Thus \(a_{\tau,m}\) measures the average amount of the original
\(F_0\)-versus-\(F_1\) ranking signal that remains visible during the
first \(m\) post-change steps.

Consider the restart-mixture process
\(\widehat M_t^{R,J} := \sum_{s=0}^{\infty} \pi_s \widehat M_t^{(s,J)}\)
from Section~\ref{subsec:restart-KT}, and define
\(\widehat D_R := \inf \left\{ t\geq1: \widehat M_t^{R,J}>B_t \right\},\)
where \(B_t>0\) is deterministic.

Since all mixture components are nonnegative,
\(\widehat M_{\tau+m}^{R,J} \geq \pi_\tau \widehat M_{\tau+m}^{(\tau,J)}.\)
Applying Lemma~\ref{lem:KT-pathwise-evidence} to the component restarted
at the true change point therefore gives
$$
    \log
    \widehat M_{\tau+m}^{R,J}
    \geq
    2m
    \left(
        \left|
            \overline p_{\tau,m}
            -
            \frac12
        \right|
        -
        \frac1{2J}
    \right)_+^2
-
    \frac{J-1}{2}\log(m+1)
    -
    C_J
    -
    \log\frac1{\pi_\tau}.
$$

The remaining ingredient is concentration of
\(\overline p_{\tau,m}\).  A bounded-difference argument applied directly
to the independent scores and auxiliary randomizers gives, for every
\(x>0\),
$$
    \mathbb P_\tau
    \left\{
        \varepsilon
        \left(
            \overline p_{\tau,m}
            -
            \mathbb E_\tau\overline p_{\tau,m}
        \right)
        \leq
        -x
    \right\}
    \leq
    \exp
    \left(
        -\frac{2mx^2}{7}
    \right).
$$

Define
$$
    \Lambda_{\tau,m,J}^{R}
    :=
    \log B_{\tau+m}
    +
    \frac{J-1}{2}\log(m+1)
    +
    C_J
    +
    \log\frac1{\pi_\tau}.
$$

\begin{theorem}[Finite-sample power of the restart-mixture KT process]
\label{thm:any-horizon-restarted-KT-power}

For every \(m\geq1\),
$$
    \mathbb P_\tau 
    \left( 
        \widehat D_R>\tau+m 
    \right) 
    \leq 
    \exp 
    \left[ 
        -\frac{2m}{7} 
        \left\{ 
            \Delta a_{\tau,m} 
            - 
            \frac1{2J} 
            - 
            \sqrt{ 
                \frac{ 
                    [\Lambda_{\tau,m,J}^{R}]_+ 
                }{ 
                    2m 
                } 
            } 
        \right\}_+^2 
    \right]. 
$$

Consequently, for every post-change monitoring horizon \(n\geq1\),
\begin{equation} 
\begin{aligned} 
    \mathbb P_\tau 
    \left( 
        \widehat D_R>\tau+n 
    \right) 
    \leq 
    \inf_{1\leq m\leq n} 
    \exp 
    \left[ 
        -\frac{2m}{7} 
        \left\{ 
            \Delta a_{\tau,m} 
            - 
            \frac1{2J} 
            - 
            \sqrt{ 
                \frac{ 
                    [\Lambda_{\tau,m,J}^{R}]_+ 
                }{ 
                    2m 
                } 
            } 
        \right\}_+^2 
    \right]. 
\end{aligned} 
\label{eq:restart-KT-any-horizon-power} 
\end{equation}
\end{theorem}

The four quantities appearing in
\eqref{eq:restart-KT-any-horizon-power} have distinct roles.  The term
$
    \Delta a_{\tau,m}
$
is the post-change rank signal,
\(1/(2J)\) is the histogram discretization cost,
$
    \frac{J-1}{2}\log(m+1)+C_J
$
is the KT learning regret, and
$
    \log\frac1{\pi_\tau}
$
is the price paid for not knowing the true change point in advance.

The finite-horizon weighting scheme of
Theorem~\ref{thm:uniform-restart-minimax},
\(\pi_s=\frac1T, \qquad s=0,\ldots,T-1,\)
makes the last term equal to
$
    \log T.
$
Hence adaptation to the unknown change point incurs only a logarithmic
penalty.

The finite-sample result also gives a simple consistency statement.
Suppose that
\(\frac{m}{\tau} \to c\in(0,\infty).\)
Then
\(a_{\tau,m} \longrightarrow \frac{\log(1+c)}{c}.\)
Therefore, if
\(\Delta \frac{\log(1+c)}{c} > \frac1{2J},\)
and
\(\log B_{\tau+m} + \log\frac1{\pi_{\tau}} = o(m),\)
then
\(\mathbb P_{\tau} \left( \widehat D_R \leq \tau+m \right) \longrightarrow1.\)
It is worth noting that the behavior of \(a_{\tau,m}\) becomes
qualitatively different when the post-change segment is much longer than
the pre-change history.  If
$
    \frac{m}{\tau}\to\infty,
$
then
$
    a_{\tau,m}
    \to0.
$
At first sight, this may appear to suggest that a longer post-change
monitoring period is detrimental to power.  This is not the correct
interpretation.  The decay of \(a_{\tau,m}\) reflects a structural
dilution of the conformal rank signal.  Indeed, a late post-change observation is ranked increasingly against
previous observations drawn from the same distribution \(F_1\).
Consequently, its conformal \(p\)-value distribution gradually moves
back toward the uniform null behavior, and averaging over a very long
post-change segment dilutes the initial \(F_0\)-versus-\(F_1\) signal.

Importantly, this does not imply that the actual detection probability
decreases with the monitoring horizon.  Since
\(\{\widehat D_R\leq\tau+m\} \subseteq \{\widehat D_R\leq\tau+n\}, \qquad m\leq n,\)
the detection probability is nondecreasing in \(n\).  This is also why
the any-horizon bound optimizes over all \(1\leq m\leq n\): even when
\(n/\tau\) is large and \(a_{\tau,n}\) is small, detection may already
have occurred at an earlier time \(m\) for which the post-change rank
signal is still strong.  Thus the decay of \(a_{\tau,m}\) describes the
weakening of the terminal-time average-\(p\) signal, rather than a loss
of power caused simply by monitoring for a longer period.

\subsubsection{Unrestarted KT process}
\label{subsubsec:unrestart-KT-power}

We now turn to the unrestarted KT process.  The same pathwise argument
applies, but the histogram is fitted to the entire history rather than
only to the observations following the change.  Define
\(\overline p_{\tau,m}^{\,\mathrm{all}} := \frac1{\tau+m} \sum_{t=1}^{\tau+m} p_t\)
and
\begin{equation}
    b_{\tau,m}
    :=
    \frac1{\tau+m}
    \sum_{k=1}^m
    \frac{\tau}{\tau+k}.
    \label{eq:b-tau-m}
\end{equation}

Since the pre-change conformal \(p\)-values have mean \(1/2\),
$$
    \varepsilon
    \left(
        \mathbb E_\tau
        \overline p_{\tau,m}^{\,\mathrm{all}}
        -
        \frac12
    \right)
    =
    \Delta b_{\tau,m}.
$$

Let \(\widehat M_t^{(J)}\) denote the unrestarted KT process from
Section~\ref{subsec:KT-histogram}, and define
\[
    \widehat D_0
    :=
    \inf
    \left\{
        t\geq1:
        \widehat M_t^{(J)}>B_t
    \right\}.
\]
Applying Lemma~\ref{lem:KT-pathwise-evidence} to all
\(\tau+m\) conformal \(p\)-values gives
$$\log
    \widehat M_{\tau+m}^{(J)}
    \geq
    2(\tau+m)
    \left(
        \left|
            \overline p_{\tau,m}^{\,\mathrm{all}}
            -
            \frac12
        \right|
        -
        \frac1{2J}
    \right)_+^2
    -
    \frac{J-1}{2}
    \log(\tau+m+1)
    -
    C_J.$$
The full-history average satisfies
$$
    \mathbb P_\tau
    \left\{
        \varepsilon
        \left(
            \overline p_{\tau,m}^{\,\mathrm{all}}
            -
            \mathbb E_\tau
            \overline p_{\tau,m}^{\,\mathrm{all}}
        \right)
        \leq
        -x
    \right\}
    \leq
    \exp
    \left(
        -\frac{2(\tau+m)x^2}{7}
    \right),
$$
for every \(x>0\).

Define
$$
    \Lambda_{\tau,m,J}^{(0)}
    :=
    \log B_{\tau+m}
    +
    \frac{J-1}{2}
    \log(\tau+m+1)
    +
    C_J.
$$

\begin{theorem}[Finite-sample power of the unrestarted KT histogram process]
\label{thm:unrestarted-KT-histogram-power}

For every \(m\geq1\),
$$\mathbb P_\tau 
    \left( 
        \widehat D_0>\tau+m 
    \right) 
    \leq 
    \exp 
    \left[ 
        -\frac{2(\tau+m)}{7} 
        \left\{ 
            \Delta b_{\tau,m} 
            - 
            \frac1{2J} 
            - 
            \sqrt{ 
                \frac{ 
                    [\Lambda_{\tau,m,J}^{(0)}]_+ 
                }{ 
                    2(\tau+m) 
                } 
            } 
        \right\}_+^2 
    \right]. $$

Consequently, for every post-change monitoring horizon \(n\geq1\),
$$\mathbb P_\tau 
    \left( 
        \widehat D_0>\tau+n 
    \right) 
    \leq 
    \inf_{1\leq m\leq n} 
    \exp 
    \left[ 
        -\frac{2(\tau+m)}{7} 
        \left\{ 
            \Delta b_{\tau,m} 
            - 
            \frac1{2J} 
            - 
            \sqrt{ 
                \frac{ 
                    [\Lambda_{\tau,m,J}^{(0)}]_+ 
                }{ 
                    2(\tau+m) 
                } 
            } 
        \right\}_+^2 
    \right].$$
\end{theorem}

If
$
    \frac{m}{\tau}
    \to
    c\in(0,\infty),
$
then
\(b_{\tau,m} \longrightarrow \frac{\log(1+c)}{1+c}.\)
Hence, if
\(\Delta \frac{\log(1+c)}{1+c} > \frac1{2J}\)
and
\(\log B_{\tau+m} = o(\tau+m),\)
then
\(\mathbb P_{\tau} \left( \widehat D_0 \leq \tau+m \right) \longrightarrow1.\)

\subsubsection{Restart versus full-history learning}
\label{subsubsec:restart-vs-unrestart-power}

The two power bounds differ primarily through their deterministic signal
factors.  From
\eqref{eq:a-tau-m} and \eqref{eq:b-tau-m},
\(b_{\tau,m} = \frac{m}{\tau+m} a_{\tau,m}.\)
Thus the full-history estimator retains only the fraction
\(\frac{m}{\tau+m}\)
of the signed average-\(p\) signal available to the component restarted
at the true change point.  This identity gives a direct finite-sample
form of the pre-change contamination effect described in
Section~\ref{subsec:restart-KT}.

The contrast is particularly strong shortly after a late change.  If
$
    m=o(\tau),
$
then
$
    a_{\tau,m}=1+o(1),
$
whereas
$
    b_{\tau,m}
    =
    \frac{m}{\tau}
    \{1+o(1)\}.
$
Hence the true-restart component retains essentially the full AUC signal
\(\Delta\), while the signal available to the unrestarted histogram is
attenuated by the small factor \(m/\tau\).  In particular, for fixed
\(J\), the restarted procedure can satisfy
$
    \Delta a_{\tau,m}>\frac1{2J}
$
at a time when the unrestarted signal
$
    \Delta b_{\tau,m}
$
is still below the histogram resolution level.
If instead
$
    \frac{m}{\tau}\to c\in(0,\infty),
$
then
$
    \frac{b_{\tau,m}}{a_{\tau,m}}
    \to
    \frac{c}{1+c}.
$
The restart advantage thus persists on proportional post-change
horizons, subject to the additional logarithmic aggregation penalty
\(\log(1/\pi_\tau)\).

The bounds also reveal the familiar resolution--estimation tradeoff in
the number of histogram bins.  Increasing \(J\) decreases the
discretization term
$
    \frac1{2J},
$
but increases the KT regret
$
    \frac{J-1}{2}\log N.
$
A finer histogram can therefore detect weaker departures from uniformity
once sufficiently many observations have accumulated, whereas a coarser
histogram pays a smaller online learning penalty at short horizons.

Finally, the rejection boundary enters both power bounds only through
\(\log B_{\tau+m}\).  Consequently, replacing the Ville boundary
\(1/\alpha\) by a smaller reference-null calibrated boundary decreases
the evidence requirement in both the restarted and unrestarted
procedures without altering their betting strategies.

\subsection{High-probability detection-delay gain}
\label{sec:delay-gain}

The preceding power analysis quantifies the benefit of a smaller
rejection boundary through the probability of detection before a
prescribed horizon. We now study a complementary question. Suppose that
the calibrated rule has already rejected while the Ville rule has not.
How much additional time can the same KT process require before reaching
the Ville boundary?

As in the preceding analysis, we condition on the independent
reference-null sample used to construct the calibrated boundary. The
boundary is therefore fixed and independent of the test trajectory.

For a given KT process, let $L_t$ denote its log wealth and let
$D_{\mathrm C}$ and $D_{\mathrm V}$ be the first strict crossings of
the calibrated boundary $B_t^{\mathrm C}$ and the Ville boundary
$1/\alpha$, respectively. We consider the event
\(D_{\mathrm C}<D_{\mathrm V},\)
which was shown in Theorem~\ref{thm:boundary-reduction} to occur with
high probability under an effective reference-null calibration.
Writing
\(d:=D_{\mathrm C},\)
define
\(h_d := \log\frac1\alpha-L_d.\)
Thus $h_d$ is exactly the additional amount of log wealth required at
the calibrated crossing time to reach the Ville boundary. If
\(B_d^{\mathrm C} = \frac{\rho_d}{\alpha}, \qquad 0<\rho_d<1,\)
and
\(O_{\mathrm C} := L_d-\log B_d^{\mathrm C}\)
denotes the overshoot above the calibrated boundary, then
\(h_d = \log\frac1{\rho_d}-O_{\mathrm C}.\)
Hence $\log(1/\rho_d)$ is the log-scale advantage created by calibration,
after subtracting the part already consumed by the overshoot.

The delay problem is therefore reduced to controlling the future
increase of the KT log wealth. We first develop the argument for the
unrestarted process, for which the main mechanism can be seen most
transparently, and then extend it componentwise to the restart mixture.

\subsubsection{Unrestarted KT process}
\label{sec:unrestarted-delay}

We first consider the unrestarted KT histogram process, for which the
mechanism behind the delay gain is most transparent. Let
$\widehat M_t^{(J)}$ be the process defined in Section~3.2 and write
\(L_t:=\log \widehat M_t^{(J)}.\)
Let
\[
    D_{\mathrm C}
    :=
    \inf\left\{
        t\geq\tau+1:
        \widehat M_t^{(J)}>B_t^{\mathrm C}
    \right\},
    \qquad
    D_{\mathrm V}
    :=
    \inf\left\{
        t\geq\tau+1:
        \widehat M_t^{(J)}>\frac1\alpha
    \right\},
\]
where both crossings are defined by strict exceedance. We focus on the
event
\(D_{\mathrm C}<D_{\mathrm V}\)
and write
\(d:=D_{\mathrm C}.\)
At time $d$, the calibrated rule has already rejected whereas the Ville
rule still requires
\(h_d := \log\frac1\alpha-L_d \geq0\)
additional units of log wealth.

If
\(B_d^{\mathrm C} = \frac{\rho_d}{\alpha}, \qquad 0<\rho_d<1,\)
and
\(O_{\mathrm C} := L_d-\log B_d^{\mathrm C}\)
is the overshoot above the calibrated boundary, then
\(h_d = \log\frac1{\rho_d}-O_{\mathrm C}.\)
Thus the reduction of the rejection boundary creates a log-evidence
advantage $\log(1/\rho_d)$, of which $O_{\mathrm C}$ has already been
used up at the calibrated crossing.

The problem is therefore to determine how quickly the KT process can
accumulate the remaining amount $h_d$ after time $d$. To this end, let
\[
    \widehat{\boldsymbol\theta}_t
    :=
    \left(
        \frac{N_{t,1}}{t},
        \ldots,
        \frac{N_{t,J}}{t}
    \right),
    \qquad
    \boldsymbol u_J
    :=
    \left(
        \frac1J,\ldots,\frac1J
    \right),
\]
and write the KT log wealth as
\begin{equation*}
\label{eq:unrestart-KT-regret-representation}
    L_t
    =
    tD_{\mathrm{KL}}
    \left(
        \widehat{\boldsymbol\theta}_t
        \,\middle\|\,
        \boldsymbol u_J
    \right)
    -
    R_t,
\end{equation*}
where $R_t$ denotes the KT regret.

For $m\geq1$, let
\[
    \widehat{\boldsymbol r}_{d,m}
    :=
    \frac1m
    \sum_{i=1}^m
    \left(
        \mathbf 1\{p_{d+i}\in I_1\},
        \ldots,
        \mathbf 1\{p_{d+i}\in I_J\}
    \right)
\]
be the empirical histogram of the next $m$ conformal $p$-values. Since
\[
    \widehat{\boldsymbol\theta}_{d+m}
    =
    \frac{d}{d+m}\widehat{\boldsymbol\theta}_d
    +
    \frac{m}{d+m}\widehat{\boldsymbol r}_{d,m},
\]
the weighted KL identity gives
\begin{equation}
\label{eq:unrestart-core-decomposition}
\begin{split}
    L_{d+m}-L_d
    ={}&
    mD_{\mathrm{KL}}
    \left(
        \widehat{\boldsymbol r}_{d,m}
        \,\middle\|\,
        \boldsymbol u_J
    \right)
    -
    \mathcal J_{d,m}
    -
    (R_{d+m}-R_d),
\end{split}
\end{equation}
where
\begin{equation}
\label{eq:unrestart-Jensen-gap}
\begin{split}
    \mathcal J_{d,m}
    :={}&
    dD_{\mathrm{KL}}
    \left(
        \widehat{\boldsymbol\theta}_d
        \,\middle\|\,
        \widehat{\boldsymbol\theta}_{d+m}
    \right)
    \\
    &+
    mD_{\mathrm{KL}}
    \left(
        \widehat{\boldsymbol r}_{d,m}
        \,\middle\|\,
        \widehat{\boldsymbol\theta}_{d+m}
    \right).
\end{split}
\end{equation}

Equation~\eqref{eq:unrestart-core-decomposition} is the key structural
identity. The first term is the new evidence supplied by the future
conformal $p$-values. The second is a nonnegative Jensen gap measuring
the incompatibility between the histogram already accumulated at time
$d$ and the future histogram. The last term is the additional KT
learning regret. We now control these three terms in turn.

\paragraph{Future evidence.}

The first term in the decomposition
\eqref{eq:unrestart-core-decomposition} measures the amount of new
evidence that can be generated by future conformal $p$-values. A
distinctive feature of conformal monitoring is that the post-change
$p$-values do not follow a fixed alternative distribution. Indeed, as
more observations from $F_1$ enter the reference sample, a future
observation $X_t\sim F_1$ is ranked against an empirical distribution
that becomes increasingly close to $F_1$. Consequently, the resulting
conformal ranks gradually return towards their null uniform
distribution.

We first quantify this phenomenon by controlling the deviation of the
future empirical histogram from the uniform histogram. Define
\begin{equation*}
\label{eq:unrestart-bd}
    b_d
    :=
    \sup_x
    \left|
        \frac1d
        \sum_{\ell=1}^d
        \mathbf 1\{X_\ell\leq x\}
        -
        F_1(x)
    \right|.
\end{equation*}
The quantity $b_d$ measures the discrepancy between the empirical
distribution used for conformal ranking at time $d$ and the future
distribution $F_1$. When $b_d$ is small, future observations are ranked
against a reference distribution close to their own generating
distribution, and hence their conformal $p$-values are close to uniform.

The following lemma provides a simultaneous high-probability control of
the future histogram.

\begin{lemma}[Uniform deviation of the future conformal histogram]
\label{lem:future-conformal-histogram}

For $k\geq1$ and $\eta\in(0,1)$, define
\[
    \lambda_k(\eta)
    :=
    \sqrt{
        \frac1{2k}
        \log\frac{2\pi^2k^2}{3\eta}
    },
    \qquad
    \lambda_0(\eta):=0 .
\]
For $b\geq0$, define
\begin{equation}
\label{eq:unrestart-adm}
\begin{split}
    a_{d,m}(b;\eta)
    :={}&
    \frac2m
    \sum_{i=1}^m
    \left[
        \frac{d}{d+i-1}b
        +
        \frac{i-1}{d+i-1}
        \lambda_{i-1}(\eta)
        +
        \frac1{d+i}
    \right]
    \\
    &+
    \sqrt{
        \frac1{2m}
        \log
        \frac{2J\pi^2m^2}{3\eta}
    } .
\end{split}
\end{equation}
Then, conditionally on $\mathcal{H}_d:=\sigma(X_1,\ldots,X_d, U_1,\ldots, U_d)$, with
probability at least $1-\eta$,
\begin{equation}
\label{eq:unrestart-future-histogram}
    \left\|
        \widehat{\boldsymbol r}_{d,m}
        -
        \boldsymbol u_J
    \right\|_\infty
    \leq
    a_{d,m}(b_d;\eta)
\end{equation}
simultaneously for all $m\geq1$.

\end{lemma}

The bound follows from two sources of randomness. First, the empirical
distribution used in conformal ranking differs from $F_1$ because the
reference sample still contains pre-change observations and because the
post-change observations in the reference sample are empirical rather
than exactly distributed according to $F_1$. The corresponding deviation
is bounded by
\(\frac{d}{d+i-1}b_d + \frac{i-1}{d+i-1} \lambda_{i-1}(\eta),\)
where the first term represents the remaining contribution of the
pre-change distribution and the second term follows from a time-uniform
DKW bound for the post-change observations.

Second, the finite-sample definition of conformal ranks introduces an
additional correction of order
\(\frac1{d+i}.\)
Averaging these deviations over the first $m$ future observations gives
the first term in \eqref{eq:unrestart-adm}. Finally, conditional on the
individual conformal $p$-values, the multinomial fluctuation of the
future histogram over the $J$ bins is controlled by a union bound and
Hoeffding's inequality, yielding the second term in
\eqref{eq:unrestart-adm}.

The lemma shows that $a_{d,m}(b_d;\eta)$ is a high-probability upper
bound for the distance between the future empirical histogram and the
uniform histogram. We now translate this histogram deviation into a
bound on the largest KL evidence that can be generated by future
observations.

Define
$$
    \Phi_J(a)
    :=
    \sup_{\substack{
        \boldsymbol q\in\Delta_J\\
        \|\boldsymbol q-\boldsymbol u_J\|_\infty\leq a
    }}
    D_{\mathrm{KL}}
    \left(
        \boldsymbol q
        \,\middle\|\,
        \boldsymbol u_J
    \right).
$$
Thus, $\Phi_J(a)$ represents the largest possible histogram evidence
against uniformity among all histograms contained in an
$\ell_\infty$-neighbourhood of radius $a$ around $\boldsymbol u_J$.
In particular,
$$
    \Phi_J(a)
    \leq
    \min\{\log J,J^2a^2\}.
$$
Combining the lemma with the definition of $\Phi_J$, we obtain, with
probability at least $1-\eta$,
\begin{equation}
\label{eq:unrestart-future-KL}
    D_{\mathrm{KL}}
    \left(
        \widehat{\boldsymbol r}_{d,m}
        \,\middle\|\,
        \boldsymbol u_J
    \right)
    \leq
    \Phi_J
    \left(
        a_{d,m}(b_d;\eta)
    \right)
\end{equation}
simultaneously for all $m\geq1$.

\paragraph{Contraction of the accumulated evidence.}

We next lower bound the Jensen gap in
\eqref{eq:unrestart-core-decomposition}. Define
\begin{equation}
\label{eq:unrestart-Vd}
    V_d
    :=
    \operatorname{TV}
    \left(
        \widehat{\boldsymbol\theta}_d,
        \boldsymbol u_J
    \right).
\end{equation}
The quantity $V_d$ measures the amount of histogram non-uniformity
already accumulated by the KT learner when the calibrated rule
rejects.

Pinsker's inequality applied to
\eqref{eq:unrestart-Jensen-gap} gives
\[
    \mathcal J_{d,m}
    \geq
    \frac{2dm}{d+m}
    \operatorname{TV}
    \left(
        \widehat{\boldsymbol\theta}_d,
        \widehat{\boldsymbol r}_{d,m}
    \right)^2.
\]
On the event \eqref{eq:unrestart-future-histogram},
\[
    \operatorname{TV}
    \left(
        \widehat{\boldsymbol r}_{d,m},
        \boldsymbol u_J
    \right)
    \leq
    \frac J2a_{d,m}(b_d;\eta),
\]
and therefore
\begin{equation}
\label{eq:unrestart-contraction}
    \mathcal J_{d,m}
    \geq
    \frac{2dm}{d+m}
    \left[
        V_d
        -
        \frac J2a_{d,m}(b_d;\eta)
    \right]_+^2.
\end{equation}

The meaning of this term is important. If the histogram at the
calibrated crossing is appreciably non-uniform but the future histogram
is already close to uniform, then the two empirical distributions are
necessarily separated. Pooling the new observations with the old
histogram therefore contracts the empirical KL divergence. Thus,
future observations may dilute rather than reinforce the evidence that
triggered calibrated detection.

\paragraph{KT regret.}

The last term in
\eqref{eq:unrestart-core-decomposition} is the increment of the KT
regret. Unlike the first two terms, this term does not describe the
distributional behavior of future conformal $p$-values. Instead, it
characterizes the additional coding cost incurred by using the KT
predictive distribution instead of the best fixed histogram model.

To control this term uniformly over the post-detection period, we first
construct a high-probability event ensuring that the histogram counts
from the pre-change uniform observations are sufficiently regular.

For $\eta_0\in(0,1)$, define
\(c_{\tau,J}(\eta_0) := \frac1J - \sqrt{ \frac1{2\tau} \log\frac{J}{\eta_0} } .\)
The condition
\(c_{\tau,J}(\eta_0)>0\)
ensures that the deviation bound below is non-vacuous. Define the event
\begin{equation}
\label{eq:unrestart-A-tau}
    \mathcal A_\tau(\eta_0)
    :=
    \left\{
        \min_{1\leq j\leq J}
        \frac{N_{\tau,j}}{\tau}
        \geq
        c_{\tau,J}(\eta_0)
    \right\}.
\end{equation}
Since the pre-change conformal $p$-values are i.i.d.
$\operatorname{Unif}(0,1)$, the corresponding histogram counts follow
a multinomial distribution. Therefore, a standard concentration
argument gives the following result.

\begin{lemma}[Uniform occupancy of the pre-change histogram]
\label{lem:KT-occupancy}

For any $\eta_0\in(0,1)$ satisfying
$c_{\tau,J}(\eta_0)>0$,
\(\mathbb P_\tau \left( \mathcal A_\tau(\eta_0) \right) \geq 1-\eta_0 .\)

\end{lemma}

The event $\mathcal A_\tau(\eta_0)$ guarantees that every histogram bin
contains a sufficiently large number of observations before the change.
This lower bound prevents the KT predictive probabilities from becoming
too small and allows the finite-sample regret expansion to hold
uniformly over all subsequent times.

To state the resulting regret bound, define
\begin{equation}
\label{eq:unrestart-K-tau-J}
\begin{split}
    K_{\tau,J}(\eta_0)
    :={}&
    \frac{J(J-1)}4
    +
    \frac{J^2(J+1)}{16}
    +
    \frac1{12}
    +
    \frac{5J}{
        24c_{\tau,J}(\eta_0)
    } .
\end{split}
\end{equation}
The quantity $K_{\tau,J}(\eta_0)$ collects the finite-sample constants
arising from the Stirling approximation and the lower bound on the
initial bin counts. Define further
\(\varepsilon_{\tau,J}(\eta_0) := \frac{2K_{\tau,J}(\eta_0)}{\tau}.\)
For fixed $J$, this quantity is of order $O(\tau^{-1})$.

\begin{lemma}[Uniform KT-regret increment bound]
\label{lem:KT-regret-increment}

On the event $\mathcal A_\tau(\eta_0)$,
the KT regret increment satisfies
\begin{equation}
\label{eq:unrestart-regret-increment}
    R_{d+m}-R_d
    \geq
    \frac{J-1}{2}
    \log\left(
        1+\frac md
    \right)
    -
    \varepsilon_{\tau,J}(\eta_0)
\end{equation}
simultaneously for all $m\geq1$.

\end{lemma}

The lemma shows that the regret accumulated after the calibrated
crossing grows only logarithmically with the relative increase of the
sample size. In particular,
\(R_{d+m}-R_d = O_J \left( \log(1+m/d) \right)\)
up to a finite-sample remainder of order $O(\tau^{-1})$.

Hence, compared with the two distribution-dependent terms in
\eqref{eq:unrestart-core-decomposition}, the KT-regret contribution is
a lower-order correction in the post-detection growth analysis. It
enters only through the explicit penalty
$\varepsilon_{\tau,J}(\eta_0)$ and the logarithmic regret increment
above.

Combining
\eqref{eq:unrestart-future-KL},
\eqref{eq:unrestart-contraction}, and
\eqref{eq:unrestart-regret-increment}, define
\begin{equation}
\label{eq:unrestart-U}
    U_{d,m}
    (b,V;\eta_0,\eta)
    :={}
    m
    \Phi_J
    \left(
        a_{d,m}(b;\eta)
    \right)-
    \frac{2dm}{d+m}
    \left[
        V
        -
        \frac J2a_{d,m}(b;\eta)
    \right]_+^2-
    \frac{J-1}{2}
    \log\left(1+\frac md\right)
    +
    \varepsilon_{\tau,J}(\eta_0).
\end{equation}
Thus $U_{d,m}$ is an upper confidence bound on the total increase of the
KT log wealth during the first $m$ observations after calibrated
detection.

\begin{theorem}[High-probability delay gain for the unrestarted KT process]
\label{thm:unrestart-delay-gain}

On $\mathcal A_\tau(\eta_0)$, conditionally on $\mathcal H_d$, with
probability at least $1-\eta$,
\(L_{d+m}-L_d \leq U_{d,m} (b_d,V_d;\eta_0,\eta)\)
simultaneously for all $m\geq1$.

Define
\begin{equation}
\label{eq:unrestart-rd}
    r_d(\eta_0,\eta)
    :=
    \sup
    \left\{
        r\in\mathbb N_0:
        \max_{1\leq m\leq r}
        U_{d,m}
        (b_d,V_d;\eta_0,\eta)
        <
        h_d
    \right\}.
\end{equation}
Then, on $\mathcal A_\tau(\eta_0)$,
\begin{equation}
\label{eq:unrestart-delay-result}
    \mathbb P_\tau
    \left(
        D_{\mathrm V}-D_{\mathrm C}
        >
        r_d(\eta_0,\eta)
        \,\middle|\,
        \mathcal H_d
    \right)
    \geq
    1-\eta.
\end{equation}
\end{theorem}

Theorem~\ref{thm:unrestart-delay-gain} converts the boundary reduction
into a delay guarantee in a direct way. At time $d$, the Ville rule
still has to close the log-evidence gap $h_d$. For each future horizon
$m$, $U_{d,m}$ bounds how much of this gap can have been closed by that
time. Hence $r_d$ is the largest window over which the upper bound on
future wealth growth remains strictly below the remaining Ville gap.
On the stated high-probability event, the Ville rule must therefore use
more than $r_d$ additional observations.

For a fixed alternative, the finite-sample value of the guaranteed
delay is obtained directly from
\eqref{eq:unrestart-rd}: evaluate $U_{d,m}$ successively as $m$ grows
and locate the first time at which it reaches $h_d$. Thus the theorem
does not merely assert that the calibrated rule is earlier; it provides
a quantitative lower bound, in observations, on how much earlier it is.

A simple approximation makes the size of this gain more transparent.
Suppose that the relevant post-detection window satisfies $m=o(d)$ and
that $a_{d,m}(b_d;\eta)$ changes slowly enough over this window to be
represented by a local value $a_d$. Then
\[
    \frac{2d}{d+m}
    =
    2+o(1),
    \qquad
    \log\left(1+\frac md\right)
    =
    \frac md+o\left(\frac md\right).
\]
Ignoring the small finite-sample remainder
$\varepsilon_{\tau,J}(\eta_0)$, the growth bound is approximately
linear:
$
    U_{d,m}
    \approx
    m\gamma_d,
$
where
\(\gamma_d := \Phi_J(a_d) - 2 \left[ V_d-\frac J2a_d \right]_+^2 - \frac{J-1}{2d}.\)
When $\gamma_d>0$, equating this approximate growth with the remaining
Ville gap gives the delay scale
\begin{equation}
\label{eq:unrestart-delay-scale}
    r_d
    \approx
    \frac{
        \log(1/\rho_d)-O_{\mathrm C}
    }{
        \gamma_d
    }.
\end{equation}

Equation~\eqref{eq:unrestart-delay-scale} gives a useful interpretation
of the theorem: the delay gain is approximately the remaining
log-boundary advantage divided by the rate at which that advantage can
be eroded after calibrated detection. 

The denominator also has a clear interpretation. A smaller $a_d$ means
that the future conformal histogram is closer to uniform, which both
reduces the new-evidence term $\Phi_J(a_d)$ and strengthens the
contraction term. A larger $V_d$ means that more non-uniformity has
already been accumulated at detection, again strengthening contraction.
Both effects decrease $\gamma_d$ and therefore increase the delay gain.

Theorem~\ref{thm:unrestart-delay-gain} provides a conditional delay
guarantee depending on the realized quantities $b_d$ and $V_d$. We next
remove this dependence by replacing them with explicit finite-sample
bounds.

\begin{corollary}
\label{cor:unrestart-delay-gain}

Let
$
    \Delta_{\mathrm{KS}}
    :=
    \sup_x|F_0(x)-F_1(x)|.
$
For $\eta_1\in(0,1)$, define
$\overline b_n(\eta_1)
:=
\frac{\tau}{n}
\{\Delta_{\mathrm{KS}}
+
\sqrt{(2\tau)^{-1}\log(4/\eta_1)}\}
+
\frac{n-\tau}{n}\lambda_{n-\tau}(\eta_1)$
for $n\geq\tau+1$.
$\kappa_J
:=
J\log\Gamma(1/2)-\log\Gamma(J/2)
+(1-J)\log(2\pi)/2$
and define
$\underline V_d^{\mathrm C}
:=
\{(2J)^{-1}
[
(\log B_d^{\mathrm C}
+(J-1)\log d/2
+\kappa_J
-K_{\tau,J}(\eta_0)/\tau)/d
]_+
\}^{1/2}$.
Then, with probability at least
$1-\eta_0-\eta_1$,
\(b_d\leq \overline b_d(\eta_1), \quad V_d\geq \underline V_d^{\mathrm C}.\)
Consequently, define
$\overline U_{d,m}(\eta_0,\eta_1,\eta_2)
:=
U_{d,m}
(\overline b_d(\eta_1),\underline V_d^{\mathrm C};
\eta_0,\eta_2)$
and
$\overline r_d(\eta_0,\eta_1,\eta_2)
:=
\sup\{r\in\mathbb N_0:
\max_{1\leq m\leq r}
\overline U_{d,m}(\eta_0,\eta_1,\eta_2)
<h_d\}$.
Then
$$
    \mathbb P_\tau
    \left(
        D_{\mathrm V}-D_{\mathrm C}
        >
        \overline r_{D_{\mathrm C}}
        (\eta_0,\eta_1,\eta_2)
    \right)
    \geq
    1-\eta_0-\eta_1-\eta_2 .
$$
\end{corollary}

The corollary provides a fully explicit delay guarantee without relying
on the realized state of the KT learner at the calibrated crossing.
Compared with the conditional result, it is less adaptive but depends
only on observable calibration information, the change magnitude through
$\Delta_{\mathrm{KS}}$, and the prescribed confidence levels.

\subsubsection{Restart-mixture KT process}

We now extend the preceding delay analysis to the restart-mixture KT
process.  The underlying principle is the same as in the unrestarted
case.  At the calibrated crossing time, the Ville rule still has to
close a remaining log-evidence gap, and the delay gain is obtained by
upper bounding how quickly the same e-process can accumulate additional
log wealth.

The main difference is that the restart mixture contains components
with different starting times.  At the calibrated crossing, some
components have already accumulated a histogram over a substantial
history, whereas other components have not yet started.  Their
post-detection behavior is therefore qualitatively different and must
be treated componentwise before being aggregated at the mixture level.

Recall that
\[
    \widehat M_t^{\mathrm R,J}
    =
    \sum_{s=0}^{\infty}
    \pi_s\widehat M_t^{(s,J)},
    \sum_{s=0}^{\infty}\pi_s=1,
\]
where
$
    \widehat M_t^{(s,J)}=1,
    \qquad t\le s.
$
Write
$
    L_t^{\mathrm R}
    :=
    \log\widehat M_t^{\mathrm R,J}.
$
Let
\[
    D_{\mathrm C}^{\mathrm R}
    :=
    \inf
    \left\{
        t\ge\tau+1:
        \widehat M_t^{\mathrm R,J}>B_t^{\mathrm C}
    \right\},
\]
and
\[
    D_{\mathrm V}^{\mathrm R}
    :=
    \inf
    \left\{
        t\ge\tau+1:
        \widehat M_t^{\mathrm R,J}>\frac1\alpha
    \right\}.
\]
On the event
$\{ D_{\mathrm C}^{\mathrm R}
    <
    D_{\mathrm V}^{\mathrm R},\}
$
write$
    d:=D_{\mathrm C}^{\mathrm R},
    h_d^{\mathrm R}
    :=
    \log\frac1\alpha-L_d^{\mathrm R}.$
Thus \(h_d^{\mathrm R}\ge0\) is the additional amount of log wealth
required by the mixture at the calibrated crossing time in order to
reach the Ville boundary.

For an active component \(s<d\), let$
    n_{d,s}:=d-s,
$/
and define its empirical histogram at time \(d\) by
\[
    \widehat{\boldsymbol\theta}_d^{(s)}
    :=
    \left(
        \frac{N_{d,1}^{(s)}}{n_{d,s}},
        \ldots,
        \frac{N_{d,J}^{(s)}}{n_{d,s}}
    \right).
\]
As before, let
\(\boldsymbol u_J := \left( \frac1J,\ldots,\frac1J \right),\)
and define
\[
    V_{d,s}
    :=
    \operatorname{TV}
    \left(
        \widehat{\boldsymbol\theta}_d^{(s)},
        \boldsymbol u_J
    \right).
\]
The relative contribution of the \(s\)-th component to the mixture
wealth at time \(d\) is
\(w_{d,s} := \frac{ \pi_s\widehat M_d^{(s,J)} }{ \widehat M_d^{\mathrm R,J} }.\)
Then$
    w_{d,s}\ge0,
    \sum_{s=0}^{\infty}w_{d,s}=1.
$

The future-distribution argument from the unrestarted case continues
to apply, but a restart component created after time \(d\) observes a
suffix rather than the complete post-\(d\) sequence.  We therefore need
the corresponding deviation bound simultaneously over all future
suffixes.

For \(r\ge0\) and \(\ell\ge1\), let$
    \widehat{\boldsymbol r}_{d,r,\ell}
$
denote the empirical \(J\)-bin histogram of
$
    (p_{d+r+1},\ldots,p_{d+r+\ell}).
$
Using the same \(\lambda_k(\eta)\) as in the preceding subsection,
define
\[
    \delta_{d,i}(b;\eta)
    :=
    2
    \left[
        \frac{d}{d+i-1}b
        +
        \frac{i-1}{d+i-1}\lambda_{i-1}(\eta)
        +
        \frac1{d+i}
    \right],
\]
and
\begin{equation}
\label{eq:RKT-a-suffix}
\begin{aligned}
    a_{d,r,\ell}^{\mathrm R}(b;\eta)
    :=
    \frac1\ell
    \sum_{i=r+1}^{r+\ell}
    \delta_{d,i}(b;\eta)
    +
    \sqrt{
        \frac1{2\ell}
        \log
        \frac{
            J\pi^4(r+1)^2\ell^2
        }{
            9\eta
        }
    }.
\end{aligned}
\end{equation}
For \(r=0\), write
\(a_{d,m}^{\mathrm R}(b;\eta) := a_{d,0,m}^{\mathrm R}(b;\eta).\)
The same DKW and martingale-difference argument used above, now with
an additional union bound over the suffix starting point \(r\), gives,
conditionally on \(\mathcal H_d\), with probability at least \(1-\eta\),
\begin{equation}
\label{eq:RKT-suffix-control}
    \left\|
        \widehat{\boldsymbol r}_{d,r,\ell}
        -
        \boldsymbol u_J
    \right\|_\infty
    \le
    a_{d,r,\ell}^{\mathrm R}(b_d;\eta)
\end{equation}
simultaneously for all \(r\ge0\) and \(\ell\ge1\).  Consequently,
$
    D_{\mathrm{KL}}
    \left(
        \widehat{\boldsymbol r}_{d,r,\ell}
        \,\middle\|\,
        \boldsymbol u_J
    \right)
    \le
    \Phi_J
    \left(
        a_{d,r,\ell}^{\mathrm R}(b_d;\eta)
    \right),
$
and
$
    \operatorname{TV}
    \left(
        \widehat{\boldsymbol r}_{d,r,\ell},
        \boldsymbol u_J
    \right)
    \le
    \frac J2
    a_{d,r,\ell}^{\mathrm R}(b_d;\eta).
$

The treatment of the KT regret also differs slightly from the
unrestarted case.  There, the first \(\tau\) uniform conformal
\(p\)-values ensured a high-probability lower bound on every histogram
count.  A restart component, especially one started close to time
\(d\), need not have such uniformly positive counts.  It is therefore
more convenient here to use the exact KT regret representation.

Let
\(a_J:=\frac{J-1}{2},\)
and retain the constant \(\kappa_J\) from the preceding subsection.
To separate the deterministic sample-size correction from the
cellwise count correction in the exact KT regret, define, for \(n\ge1\),
\[
    g_J(n)
    :=
    \log\Gamma\left(n+\frac J2\right)
    -
    n\log n
    +
    n
    -
    a_J\log n
    -
    \frac12\log(2\pi),
\]
and
\[
\psi(k)
:=
\begin{cases}
\displaystyle
k\log k
-
\log\Gamma(k+1/2)
-
k
+
\frac12\log(2\pi),
& k\ge1,
\\[1ex]
\displaystyle
\frac12\log2,
& k=0.
\end{cases}
\]

Both quantities are finite-sample Stirling corrections.  In particular,
\(g_J(n) = \left( \frac{J(J-2)}8+\frac1{12} \right)\frac1n + O_J(n^{-2}),\)
so that \(g_J(n)=O_J(n^{-1})\).  Moreover, for \(k\ge1\),
refined Stirling bounds give$
    0<\psi(k)<\frac1{24k},
$
and hence \(\psi(k)=O(k^{-1})\).  At an empty cell,
$
    \psi(0)=\frac12\log2.
$
Consequently,
\(0\le\psi(k)\le\frac12\log2, \qquad k\in\mathbb N_0.\)
Thus the uniform upper bound is attained only at an empty cell; once
a cell has accumulated observations, its correction decreases at
order \(1/k\).
For an active component \(s<d\), define
\(\Psi_{d,s} := \sum_{j=1}^J \psi\left(N_{d,j}^{(s)}\right).\)
Although \(\Psi_{d,s}\) is random unconditionally, it is
\(\mathcal F_d\)-measurable.  Thus, in the conditional analysis below,
it is part of the observed state of the \(s\)-th KT learner at the
calibrated crossing.

The exact KT regret identity implies that an active component with
\(n=n_{d,s}\) observations at time \(d\) satisfies
\begin{equation*}
\label{eq:RKT-regret-old}
\begin{aligned}
    R_{n+m}^{(s)}-R_n^{(s)}
    \ge\;&
    a_J
    \log\left(1+\frac mn\right)
    +
    g_J(n+m)-g_J(n)
    -
    \Psi_{d,s}.
\end{aligned}
\end{equation*}
For a newly initialized component based on \(\ell\) observations,
\(R_\ell^{\mathrm{new}} \ge a_J\log\ell + \kappa_J + g_J(\ell).\)

We now turn to the main distinction between the different restart
components.

Consider first a component with \(s<d\).  It has already accumulated
\(n_{d,s}=d-s\) observations before calibrated detection.  Hence its
future log-wealth increment compares a pre-existing histogram
\(\widehat{\boldsymbol\theta}_d^{(s)}\) with the histogram of the new
observations.  Exactly as in the unrestarted analysis, the weighted KL
identity therefore produces a contraction term.  Combining the
future-histogram bound with \eqref{eq:RKT-regret-old}, define
\begin{equation}
\label{eq:RKT-G-old}
\begin{aligned}
G_{d,m}^{(s)}(b;\eta)
:=\;&
m
\Phi_J
\left(
    a_{d,m}^{\mathrm R}(b;\eta)
\right)-
\frac{
    2n_{d,s}m
}{
    n_{d,s}+m
}
\left[
    V_{d,s}
    -
    \frac J2
    a_{d,m}^{\mathrm R}(b;\eta)
\right]_+^2
\\
&-
a_J
\log
\left(
    1+\frac m{n_{d,s}}
\right)-
\left\{
    g_J(n_{d,s}+m)-g_J(n_{d,s})
\right\}
+
\Psi_{d,s}.
\end{aligned}
\end{equation}
Thus, on the event \eqref{eq:RKT-suffix-control},
\[
    \log
    \frac{
        \widehat M_{d+m}^{(s,J)}
    }{
        \widehat M_d^{(s,J)}
    }
    \le
    G_{d,m}^{(s)}(b_d;\eta).
\]

The situation is different for a component with
$
    s=d+r,
    0\le r<m.
$
Such a component has not yet started at time \(d\), so
$
    \widehat M_d^{(d+r,J)}=1.
$
By time \(d+m\), it has observed only the suffix
$
    p_{d+r+1},\ldots,p_{d+m},
$
whose length is
$
    \ell:=m-r.
$
There is therefore no old histogram with which the future observations
can be incompatible, and hence no contraction term.  On the other
hand, the component starts learning its histogram from scratch and
must pay the full KT learning cost for a sample of size \(\ell\).
Define
\begin{equation}
\label{eq:RKT-H-new}
H_{d,r,m}(b;\eta)
:=
\ell
\Phi_J
\left(
    a_{d,r,\ell}^{\mathrm R}(b;\eta)
\right)-
a_J\log\ell
-
\kappa_J
-
g_J(\ell),
\qquad
\ell=m-r.
\end{equation}
Then
\(\log \widehat M_{d+m}^{(d+r,J)} \le H_{d,r,m}(b_d;\eta).\)

Finally, if \(s\ge d+m\), the component remains inactive throughout
\((d,d+m]\), and hence
\(\frac{ \widehat M_{d+m}^{(s,J)} }{ \widehat M_d^{(s,J)} } =1.\)

The three cases can now be combined using the posterior weights
\(w_{d,s}\).  Since
\[
    \frac{
        \widehat M_{d+m}^{\mathrm R,J}
    }{
        \widehat M_d^{\mathrm R,J}
    }
    =
    \sum_{s=0}^{\infty}
    w_{d,s}
    \frac{
        \widehat M_{d+m}^{(s,J)}
    }{
        \widehat M_d^{(s,J)}
    },
\]
define
\begin{equation}
\label{eq:RKT-U-cond}
U_{d,m}^{\mathrm R}(b;\eta)
:=
\log\Bigg\{
\sum_{s=0}^{d-1}
w_{d,s}
\exp
\left(
    G_{d,m}^{(s)}(b;\eta)
\right)+
\sum_{r=0}^{m-1}
w_{d,d+r}
\exp
\left(
    H_{d,r,m}(b;\eta)
\right)+
\sum_{s=d+m}^{\infty}
w_{d,s}
\Bigg\}.
\end{equation}

\begin{theorem}[High-probability delay gain for the restart-mixture KT process]
\label{thm:RKT-delay-gain}
Conditionally on \(\mathcal H_d\), with probability at least
\(1-\eta\),
\(L_{d+m}^{\mathrm R} - L_d^{\mathrm R} \le U_{d,m}^{\mathrm R}(b_d;\eta)\)
simultaneously for every \(m\ge1\).

Define
\begin{equation}
\label{eq:RKT-r-cond}
    r_d^{\mathrm R}(\eta)
    :=
    \sup
    \left\{
        r\in\mathbb N_0:
        \max_{1\le m\le r}
        U_{d,m}^{\mathrm R}(b_d;\eta)
        <
        h_d^{\mathrm R}
    \right\}.
\end{equation}
Then
\begin{equation*}
\label{eq:RKT-delay-cond}
    \mathbb P_\tau
    \left(
        D_{\mathrm V}^{\mathrm R}
        -
        D_{\mathrm C}^{\mathrm R}
        >
        r_d^{\mathrm R}(\eta)
        \,\middle|\,
        \mathcal H_d
    \right)
    \ge
    1-\eta.
\end{equation*}
\end{theorem}

The mixture structure is important for interpreting the resulting
delay bound.  Unlike the unrestarted process, the post-detection
growth rate is not determined by a single KT learner.  For each
restart time \(s\), define its normalized contribution over a future
window of length \(m\) by
\[
\gamma_{d,s}(m)
:=
\begin{cases}
\displaystyle
m^{-1}G_{d,m}^{(s)}(b_d;\eta),
& s<d,
\\[1.5ex]
\displaystyle
m^{-1}H_{d,s-d,m}(b_d;\eta),
& d\le s<d+m,
\\[1.5ex]
0,
& s\ge d+m.
\end{cases}
\]
Then the mixture growth bound can be written exactly as
\begin{equation}
\label{eq:RKT-rate-mixture}
U_{d,m}^{\mathrm R}(b_d;\eta)
=
\log
\left\{
    \sum_{s=0}^{\infty}
    w_{d,s}
    \exp
    \left(
        m\gamma_{d,s}(m)
    \right)
\right\}.
\end{equation}
Thus \(U_{d,m}^{\mathrm R}\) is the log moment-generating function of
the componentwise future growth rates under the posterior restart
weights at the calibrated crossing.

Define their posterior-weighted mean
\begin{equation*}
\label{eq:RKT-mean-rate}
\overline\gamma_{d,m}^{\mathrm R}
:=
\sum_{s=0}^{\infty}
w_{d,s}\gamma_{d,s}(m)
=
\sum_{s=0}^{d-1}
w_{d,s}
\frac{
G_{d,m}^{(s)}(b_d;\eta)
}{m}+
\sum_{r=0}^{m-1}
w_{d,d+r}
\frac{
H_{d,r,m}(b_d;\eta)
}{m}.
\end{equation*}
The second sum automatically accounts for the shorter exposure of a
component restarted after \(d\): if it starts at \(d+r\), then it has
only \(m-r\) observations by time \(d+m\).

When the dispersion of the quantities
\(m\gamma_{d,s}(m)\) under the weights \(w_{d,s}\) is moderate,
a cumulant expansion of \eqref{eq:RKT-rate-mixture} gives
\begin{equation*}
\label{eq:RKT-mixture-expansion}
U_{d,m}^{\mathrm R}
\approx
m\overline\gamma_{d,m}^{\mathrm R}
+
\frac{m^2}{2}
\operatorname{Var}_{w_d}
\left(
    \gamma_{d,S}(m)
\right).
\end{equation*}
The first term is the posterior-weighted average component growth,
while the second captures heterogeneity across candidate restart
times.  In particular, the log-sum-exp aggregation grows faster than
the simple weighted average whenever the componentwise rates are
heterogeneous.

Ignoring the second-order heterogeneity term for a first-order
approximation, the Ville delay is determined by
$
    m\overline\gamma_{d,m}^{\mathrm R}
    \approx
    h_d^{\mathrm R}.
$
If the weighted rate varies slowly over the relevant window, this
reduces to
$
    r_d^{\mathrm R}
    \approx
        h_d^{\mathrm R}
/
        \overline\gamma_d^{\mathrm R}.
$
Since
\(h_d^{\mathrm R} = \log\frac1{\rho_d} - O_{\mathrm C}^{\mathrm R},\)
the corresponding scale is
\begin{equation*}
\label{eq:RKT-delay-scale-rho}
    r_d^{\mathrm R}
    \approx
    \frac{
        \log(1/\rho_d)-O_{\mathrm C}^{\mathrm R}
    }{
        \overline\gamma_d^{\mathrm R}
    }.
\end{equation*}

We finally give a less-conditional version.  As in the unrestarted
case, let
\[
    \overline b_n(\eta_1)
    :=
    \frac{\tau}{n}
    \left\{
        \Delta_{\mathrm{KS}}
        +
        \sqrt{
            \frac1{2\tau}
            \log\frac4{\eta_1}
        }
    \right\}
    +
    \frac{n-\tau}{n}
    \lambda_{n-\tau}(\eta_1),
    \qquad
    n\ge\tau+1.
\]
With probability at least \(1-\eta_1\),
\(b_n\le\overline b_n(\eta_1)\)
simultaneously for all \(n\ge\tau+1\), and hence also at the random
crossing time \(d\).

To remove the remaining component-specific quantities, use
\(V_{d,s}\ge0, \qquad \Psi_{d,s} \le \frac J2\log2.\)
For a generic active-component age \(1\le n\le d\), define
\begin{equation*}
\label{eq:RKT-Gbar}
\begin{aligned}
\overline G_{n,m}^{\mathrm R}
(b;\eta)
:=
m
\Phi_J
\left(
    a_{d,m}^{\mathrm R}(b;\eta)
\right)-
a_J
\log
\left(
    1+\frac mn
\right)-
\left\{
    g_J(n+m)-g_J(n)
\right\}
+
\frac J2\log2.
\end{aligned}
\end{equation*}
Since
\(\log\left( \sum_s w_{d,s}e^{x_s} \right) \le \sup_s x_s,\)
define
\begin{equation}
\label{eq:RKT-U-less}
\overline U_{d,m}^{\mathrm R}
(\eta_1,\eta_2)
:=
\max\Bigg\{
0,
\max_{1\le n\le d}
\overline G_{n,m}^{\mathrm R}
\left(
    \overline b_d(\eta_1);
    \eta_2
\right),
\max_{0\le r<m}
H_{d,r,m}
\left(
    \overline b_d(\eta_1);
    \eta_2
\right)
\Bigg\}.
\end{equation}
Set
\begin{equation}
\label{eq:RKT-r-less}
    \overline r_d^{\mathrm R}
    (\eta_1,\eta_2)
    :=
    \sup
    \left\{
        r\in\mathbb N_0:
        \max_{1\le m\le r}
        \overline U_{d,m}^{\mathrm R}
        (\eta_1,\eta_2)
        <
        h_d^{\mathrm R}
    \right\}.
\end{equation}

\begin{corollary}[Less-conditional delay gain for the restart mixture]
\label{cor:RKT-delay-less}
For any \(\eta_1,\eta_2\in(0,1)\),
\begin{equation}
\label{eq:RKT-delay-less}
    \mathbb P_\tau
    \left(
        D_{\mathrm V}^{\mathrm R}
        -
        D_{\mathrm C}^{\mathrm R}
        >
        \overline r_{D_{\mathrm C}^{\mathrm R}}^{\mathrm R}
        (\eta_1,\eta_2)
    \right)
    \ge
    1-\eta_1-\eta_2.
\end{equation}
\end{corollary}

The conditional and less-conditional results serve the same purposes
as in the unrestarted case.  The conditional theorem retains the
realized posterior weights, component ages, accumulated
non-uniformities, and KT states at the calibrated crossing, and is
therefore substantially more adaptive.  The corollary removes these
quantities by taking a worst-case bound over active-component ages and
new restart times.  Its price is additional conservativeness, but its
inputs no longer depend on the detailed state of the restart mixture
at detection.

\section{Simulation study}
\label{sec:simulation-study}

We conducted a simulation study to answer three questions.  First, does
reference-null calibration attain the nominal finite-horizon type-I error?
Second, when it lowers the boundary of a fixed e-process, does the reduction
translate into greater post-change power or shorter detection delay?  Third,
how do these conclusions depend on the monitoring horizon and the location of
the change?  The presentation below is organised in that order: null validity
and boundary reduction, horizon-specific efficiency gains, and sensitivity to
the change-point location.

\subsection{Design and performance measures}
\label{subsec:simulation-design}

The data-generating model is specified at the score level.  Before the change,
$X_t\sim N(0,1)$; after time $\tau$, $X_t\sim N(\mu,1)$.  Each score sequence,
together with independent tie-breaking randomisers, is converted into one
global sequential randomized conformal $p$-value sequence.  The same sequence
is supplied to all procedures: the KT histogram process (KT), its
cautiously gated version (CKT) \citep{eliades2022betting}, the
restart-of-wealth mixture with a global histogram learner (PIT)
\citep{farran2026model}, and the restart mixture whose components restart both
wealth and the KT learner (RKT).
For every procedure, the calibrated version (suffix C) and the Ville version
(suffix V) use exactly the same e-process path and differ only in the crossing
boundary.  The resulting paired comparisons therefore isolate the effect of
calibration from that of the betting strategy.

CKT, PIT and RKT modify the evidence process in different ways.  CKT retains
one global KT histogram learner and does not restart wealth.  Instead, its
monitored wealth is allowed to use the current KT betting factor only after the
unrestricted KT process has risen by a factor exceeding $\epsilon$ from its
running minimum; while this cautious gate is closed, the monitored one-step
factor is one and the CKT path remains flat.  PIT addresses an unknown change
time by mixing wealth streams that start at different candidate times, but all
components continue to use the same one-step factor estimated by one global
histogram.  It therefore restarts accumulated wealth but not the information
used by the learner.  RKT also mixes candidate restart times, but each
component has its own KT counts and uses only observations arriving after that
component starts.  It restarts both wealth and the local learner, thereby
removing pre-change histogram counts from components initiated near the true
change.  The distinction is thus cautious gating for CKT, wealth restarting
with global learning for PIT, and joint wealth-and-learner restarting for RKT.

Reference-null calibration acts at a separate layer.  It is not another
betting or restart rule: for each fixed base process it estimates an upper-tail
quantile of the finite-horizon null maximum
$Z_T=\max_{1\leq t\leq T}M_t$ and replaces the universal Ville boundary
$1/\alpha$ by that process-specific threshold.  It changes neither the
conformal $p$-values nor the gate, histogram counts, mixture weights or process
path.  Accordingly, RKT denotes an e-process construction, whereas the
proposed calibrated boundary is modular and can be applied without alteration
to KT, CKT, PIT or RKT.  Across-process differences in the figures reflect the
construction of the e-process; within-process C-minus-V differences reflect
only the proposed boundary calibration.

The histogram has $J=20$ bins and the prespecified cautiousness parameter is
$\epsilon=100$.  The main horizons are
$T\in\{250,500,1000,2000\}$, with $\tau=T/2$ unless stated otherwise, and the
signal grid is
$\mu\in\{0.00,0.25,0.50,0.75,1.00,1.25,1.50,1.75,2.00\}$.  Calibrated
boundaries are estimated independently for each process and horizon using 100
reference banks, each containing $K=499$ null paths.  The 100 banks quantify
the variability due to estimating a boundary from a finite reference sample;
they are not additional evaluation paths.  Performance is evaluated on 20,000
fresh null paths and 5,000 fresh alternative paths per scenario.  Calibration
and evaluation samples are independent, while the evaluation paths are shared
across C and V to retain the natural pairing.  The nominal level is
$\alpha=0.05$, the calibrated crossing rule is strict, and no calibrated
threshold is truncated at $1/\alpha$.  Paired uncertainty intervals use 1,000
bootstrap resamples of the paired C-minus-V outcomes.

Let $D$ denote the first boundary-crossing time.  Conditional post-change
power is $\Pr(D\leq T\mid D>\tau)$, while $\Pr(D\leq\tau)$ is treated as a
pre-change false alarm rather than a detection.  Restricted mean detection
delay (RMDD) is the mean of $D-\tau$ among paths surviving to $\tau$, with a
non-detection by $T$ assigned the value $T-\tau+1$.  RMDD is therefore measured
in observations.  A positive value of
$\operatorname{RMDD}_{\rm V}-\operatorname{RMDD}_{\rm C}$ means that the
calibrated boundary detects the change earlier.

\subsection{Finite-horizon validity and boundary reduction}
\label{subsec:simulation-validity}

\begin{figure}[H]
    \centering
    \includegraphics[width=\textwidth]{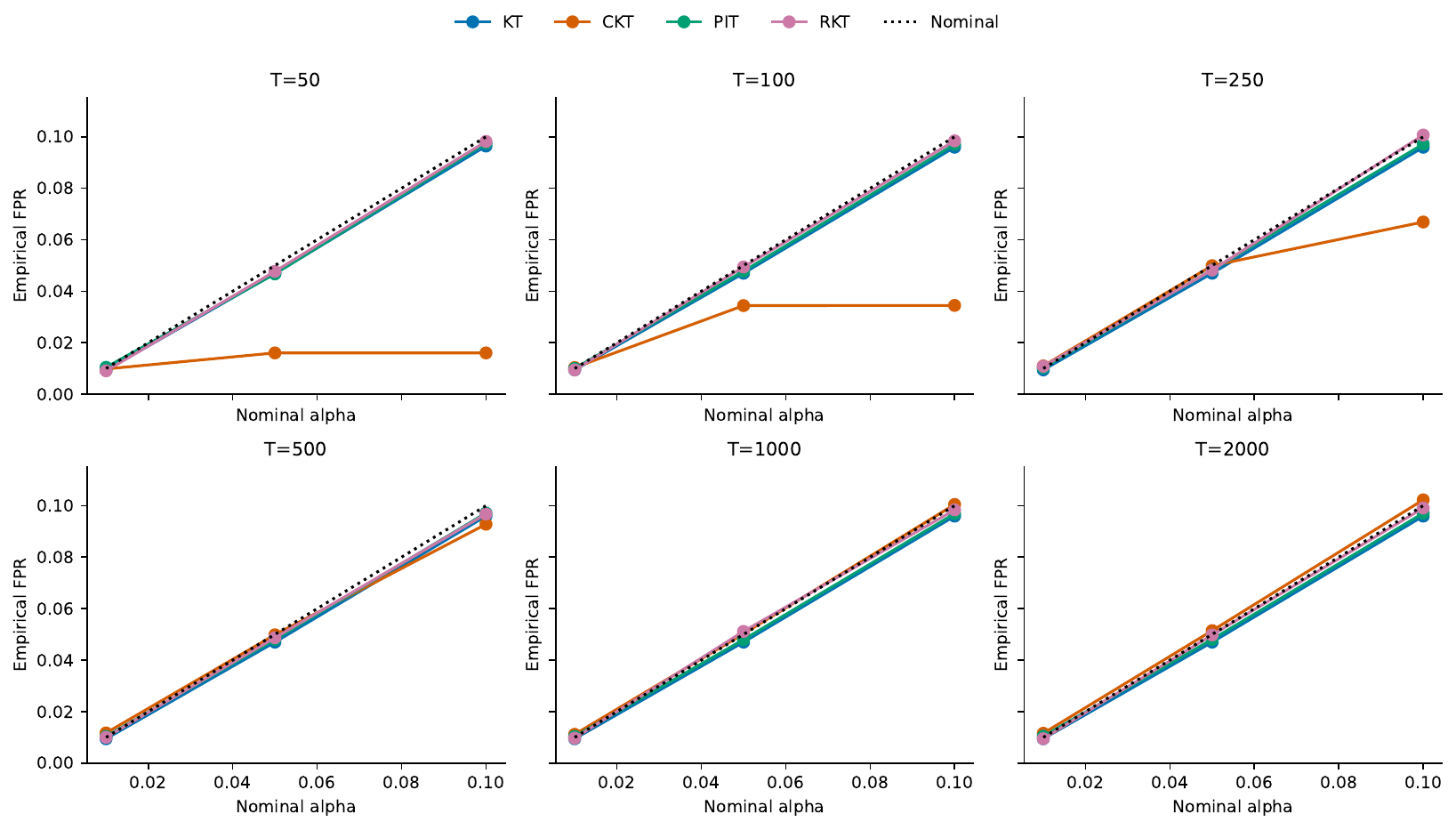}
    \caption{Empirical null crossing probabilities for the reference-null
    calibrated procedures.  Each panel corresponds to a monitoring horizon,
    and the dotted diagonal denotes equality with the nominal level.  Each
    plotted value averages the rejection probability over 100 independently
    estimated calibration boundaries, evaluated on the same 20,000 fresh null
    paths.  The separation of CKT from the diagonal at short horizons is due to
    discreteness under strict crossing, not to anti-conservative calibration.}
    \label{fig:simulation-validity}
\end{figure}

Figure~\ref{fig:simulation-validity} compares the empirical null crossing
probability with nominal levels $0.01$, $0.05$ and $0.10$ over horizons from 50
to 2,000.  The KT, PIT and RKT curves lie close to the diagonal throughout.
CKT is markedly conservative in the first three panels because its cautious
gate creates an unusually discrete short-horizon null maximum.  With
$\epsilon=100$, the gate opens only after the unrestricted KT process has
recovered more than one hundredfold from its running minimum.  Under the null,
the probability that it opens at least once is only $0.03445$, $0.06640$ and
$0.11640$ for $T=50$, $100$ and $250$, respectively; even when it opens, the
subsequent factor need not push the monitored path above one.  Consequently,
the null probabilities
\[
 \Pr\!\left(Z_T^{\rm CKT}=1\right)
 =0.98395,\ 0.96550,\ 0.93310
\]
at those horizons are very large.

This atom explains both the plateaux and their exact heights.  Whenever the
estimated upper-tail quantile is one, the strict calibrated crossing event is
$Z_T^{\rm CKT}>1$, whose probabilities are only $0.01605$, $0.03450$ and
$0.06690$ at $T=50$, $100$ and $250$.  Thus, at $T=50$, the target levels
$0.05$ and $0.10$ produce the same empirical FPR $0.01605$; at $T=100$, they
both produce approximately $0.0345$; and at $T=250$, the $0.10$ target produces
$0.0669$.  The $T=250$ target $0.05$ is attainable because its selected
quantile lies above the atom, giving empirical FPR $0.04995$.  Finite-sample
validity requires the FPR to be no greater than $\alpha$, not to equal it.
Exact equality at a level that cuts through the atom would require an
additional randomized rejection at $Z_T^{\rm CKT}=1$, which is deliberately
not used.  As $T$ increases, the gate has more opportunities to open and the
upper tail acquires finer resolution; by $T=500$ the CKT curve is close to the
diagonal, and at $T=1000$ and $T=2000$ all four procedures track the nominal
levels closely.  The first three panels therefore show benign conservativeness
caused by discreteness, rather than a calibration failure.

Table~\ref{tab:simulation-null-summary} focuses on the principal level
$\alpha=0.05$ at $T=2000$.  The calibrated false-positive rates range from
$0.0470$ to $0.0515$.  For comparison, the binomial Monte Carlo standard error
at probability $0.05$ with 20,000 paths is approximately $0.00154$, so these
departures are small in simulation terms.  The Ville boundary is conservative
for every process and especially for CKT, whose false-positive rate is only
$0.0067$.  Calibration reduces the median boundary to $65.7\%$, $14.2\%$,
$71.6\%$ and $64.4\%$ of the Ville boundary for KT, CKT, PIT and RKT,
respectively.  Because the Ville boundary is $1/\alpha=20$, the corresponding
median calibrated thresholds are approximately $13.14$, $2.84$, $14.32$ and
$12.88$.  Thus, the most conservative uncalibrated procedure, CKT, also has the
largest recoverable boundary slack.  Boundary reduction alone is not a power
claim, however; its operational value depends on how much alternative-path
mass lies between the two boundaries, as assessed below.

\begin{table}[H]
\centering
\caption{Null performance and boundary reduction at $T=2000$ and
$\alpha=0.05$.  FPR--C and FPR--V are the empirical false-positive rates for
the calibrated and Ville boundaries applied to the same 20,000 evaluation
paths.  The final column is the median, across 100 calibration banks, of the
calibrated threshold divided by the Ville threshold.}
\label{tab:simulation-null-summary}
\small
\setlength{\tabcolsep}{8pt}
\begin{tabular}{lccc}
\toprule
Process & FPR--C, $T=2000$ & FPR--V, $T=2000$ & Median threshold/Ville \\
\midrule
KT  & 0.0470 & 0.0306 & 0.657 \\
CKT & 0.0515 & 0.0067 & 0.142 \\
PIT & 0.0481 & 0.0351 & 0.716 \\
RKT & 0.0498 & 0.0320 & 0.644 \\
\bottomrule
\end{tabular}
\end{table}

\subsection{Power and detection-delay gains across horizons}
\label{subsec:simulation-horizon-gains}

\begin{figure}[H]
    \centering
    \includegraphics[width=\textwidth]{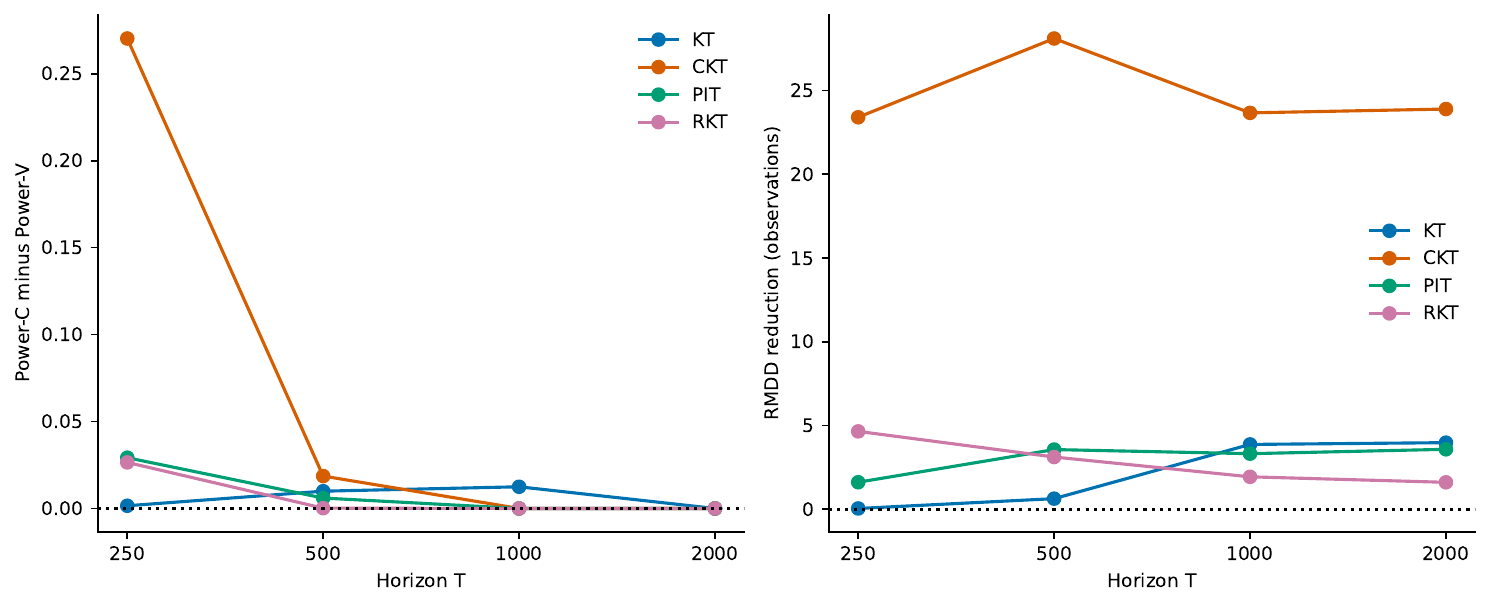}
    \caption{Gain from reference-null calibration at $\mu=1.00$ and
    $\tau=T/2$.  The left panel shows Power--C minus Power--V; the right panel
    shows $\operatorname{RMDD}_{\rm V}-\operatorname{RMDD}_{\rm C}$ in
    observations.  Positive values favour calibration.  The four horizons
    $T=250$, $500$, $1000$ and $2000$ are equally spaced categorical design
    points rather than values on a continuous or logarithmic axis.}
    \label{fig:simulation-horizon}
\end{figure}

Figure~\ref{fig:simulation-horizon} fixes $\mu=1.00$ and $\tau=T/2$ and
compares the two boundaries on the paired alternative paths.

At $T=250$,
calibration increases CKT power by $0.2703$, from $0.5059$ to $0.7762$; its
95\% paired-bootstrap interval is $(0.2594,0.2824)$.  The corresponding power
gains are smaller for PIT, RKT and KT, at $0.0292$, $0.0265$ and $0.0016$.
The reason is not that the calibration algorithm treats CKT preferentially;
the same independent-null order statistic is used for every process.  Rather,
CKT leaves substantially more unused slack in the universal Ville bound.  At
$T=250$ and $\alpha=0.05$, its median calibrated threshold is only $1.42$,
or $0.071$ of the Ville boundary 20.  The corresponding thresholds are
$13.13$, $14.26$ and $12.11$ for KT, PIT and RKT, respectively.  Hence
calibration reduces the evidence requirement roughly fourteenfold for CKT,
but by less than twofold for the other processes.

This unusually large reduction is a direct consequence of the cautious gate.
Under the null, the CKT path is usually frozen at one, so the upper null
quantile can be far below 20.  Under the alternative, a persistent departure
can open the gate and then generate sustained KT growth.  At $T=250$ many such
alternative paths have maxima between $1.42$ and 20: CKT-C detects them whereas
CKT-V does not, producing the $0.2703$ power increase.  The null maxima of KT,
PIT and RKT are more active and have heavier, more continuous upper tails, so
their calibrated thresholds remain much closer to 20; correspondingly fewer
alternative paths are reclassified by the boundary change.  Thus the power
gain depends jointly on the null boundary reduction and the amount of
alternative mass between the two thresholds, not on the threshold ratio alone.

The right panel reveals a different effect that remains visible after power
saturates.  At $T=250$, calibration reduces CKT RMDD by $23.40$ observations,
with a 95\% paired-bootstrap interval of $(22.78,24.02)$; the corresponding
reductions for KT, PIT and RKT are $0.05$, $1.61$ and $4.65$ observations.  At
$T=500$, the CKT power gain falls to $0.0186$ as both variants approach unit
power, but its RMDD reduction rises to $28.11$ observations.  At $T=1000$ and
$T=2000$, the CKT power difference is $0.0000$ to the reported precision,
whereas the RMDD improvements remain $23.67$ and $23.89$ observations.  The
other methods show the same distinction on a smaller scale: at $T=2000$, all
four power gains are $0.0000$, yet KT, PIT and RKT still detect respectively
$3.97$, $3.58$ and $1.61$ observations earlier on average.  Hence additional
power is the relevant benefit in the transition region, whereas shorter delay
is the more informative benefit once both boundaries detect almost every
surviving path.  Mechanistically, once the CKT gate opens, reaching $1.42$ (or
the corresponding calibrated threshold at a longer horizon) requires much
less post-change growth than reaching 20.  This produces a large time shift in
the crossing event even when both versions eventually cross.  Because the
other processes receive a much smaller boundary reduction, their crossing
times move less.  The persistent CKT delay gain is therefore the saturated-
power counterpart of the same process-specific Ville slack that generates its
short-horizon power gain.

\subsection{Change-point location and pre-change contamination}
\label{subsec:simulation-change-location}

\begin{figure}[H]
    \centering
    \includegraphics[width=\textwidth]{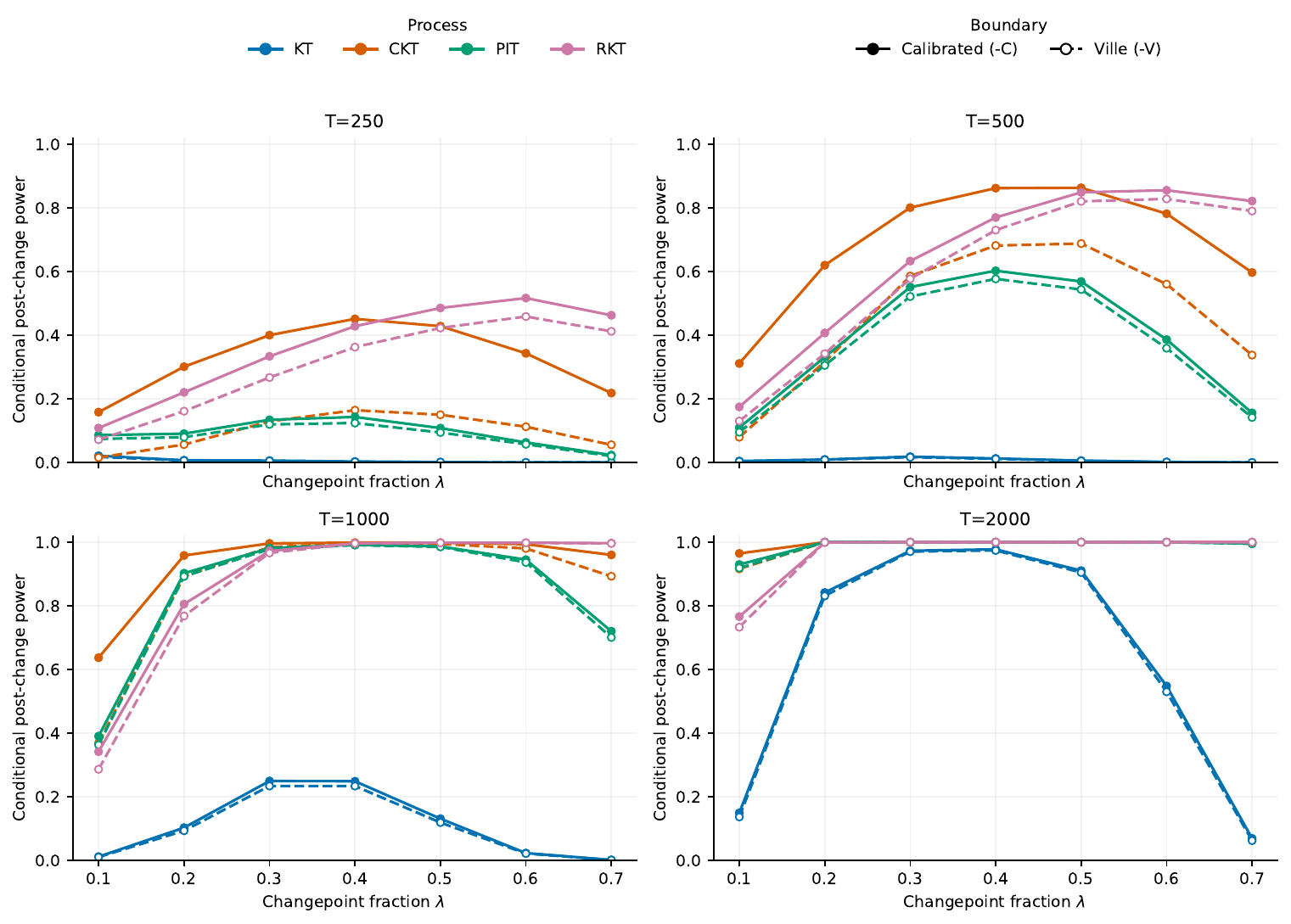}
    \caption{Conditional post-change power as a function of the change-point
    fraction $\lambda$ for a score shift $\mu=0.75$.  Panels correspond to
    $T=250$, $500$, $1000$ and $2000$; each point uses 5,000 fresh alternative
    paths, and only paths surviving to the change enter the conditional-power
    denominator.  Colours denote the e-process.  Solid curves with filled
    markers use the calibrated boundary, while dashed curves with hollow
    markers use the Ville boundary.}
    \label{fig:simulation-contamination-power}
\end{figure}

We next set $\mu=0.75$ and vary
$\tau=\operatorname{round}(\lambda T)$ over
$\lambda\in\{0.10,0.20,0.30,0.40,0.50,0.60,0.70\}$.  Figure~\ref{fig:simulation-contamination-power}
reports conditional post-change power itself, rather than only a method
difference, so both the absolute detection difficulty and the calibration gain
remain visible.  Moving the change later has two opposing effects: it provides
more pre-change observations to the learner, but leaves fewer post-change
observations before the horizon.  The resulting curves are therefore
non-monotone, particularly for KT and CKT.

At $T=250$, all methods operate in a difficult regime, but their responses are
distinct.  At $\lambda=0.40$, calibrated CKT has power $0.4513$, compared with
$0.1643$ for its Ville counterpart, while calibrated and Ville RKT attain
$0.4281$ and $0.3625$.  PIT is lower at $0.1430$ and $0.1239$, and KT is close
to zero.  By $\lambda=0.70$, RKT retains calibrated power $0.4626$, whereas CKT
falls to $0.2182$ and PIT to $0.0233$.  The restart of both wealth and the local
learner therefore helps RKT when a short post-change segment follows a long
pre-change history, although it does not establish uniform dominance.

At $T=500$, CKT reaches its highest calibrated power, approximately $0.86$, for
intermediate changes around $\lambda=0.40$--$0.50$ and still gains materially
from calibration; at $\lambda=0.40$, the calibrated and Ville powers are
$0.8621$ and $0.6817$.  RKT rises more steadily and is strongest for a late
change, reaching $0.8215$ at $\lambda=0.70$.  PIT follows the same broad shape
as RKT but at a lower level, whereas KT remains near zero.  These comparisons
show that restarting the learner, not only the wealth, is valuable in the
shorter late-change configurations.

With $T=1000$, CKT, PIT and RKT have power near one over much of the interior
range.  At $\lambda=0.30$, their calibrated powers are $0.9962$, $0.9814$ and
$0.9723$, respectively.  RKT is most robust at the late endpoint, retaining
power $0.9968$ at $\lambda=0.70$, while PIT and CKT fall to $0.7200$ and
$0.9600$.  KT peaks only around the middle and declines sharply near either
endpoint.  At $T=2000$, the three adaptive or restart-based processes are
essentially saturated after $\lambda=0.20$; even at $\lambda=0.70$, calibrated
power is $0.9998$ for CKT, $0.9960$ for PIT and $1.0000$ for RKT.  KT remains
location-sensitive: its calibrated power is about $0.1493$ at $\lambda=0.10$,
rises above $0.97$ in the interior, and falls to $0.0691$ at $\lambda=0.70$.
The C and V curves nearly coincide whenever power is saturated, while their
larger separation at $T=250$ and $T=500$ identifies the regimes in which the
lower calibrated boundary converts directly into additional detections.

Taken together, these results give a coherent finite-sample account.  The
calibrated procedures achieve approximately nominal type-I error at the main
horizons while using thresholds that are appreciably below the universal Ville
boundary.  The practical return from that slack depends on the regime: it is
seen primarily as increased power when the alternative distribution straddles
the two thresholds, and as earlier detection once power has saturated.  The
change-location experiment further shows that the magnitude of the return is
process-specific because global learning, wealth restarting and learner
restarting respond differently to pre-change history and to the remaining
post-change window.  All 34 prespecified numerical checks passed, including
checks of finite-value handling, crossing-rule consistency, figure labels,
horizon-axis formatting and trajectory provenance.

\section{Real-model experiment: online LLM watermark detection}
\label{sec:llm_experiment}

To demonstrate the practical effectiveness of the proposed reference-null
calibration method, we apply it to online large language model (LLM)
watermark detection. Although the proposed calibration procedure is
model-agnostic, LLM watermark detection provides a natural application
scenario because the detector operates sequentially on generated tokens and
the evidence accumulates over time. In this experiment, we follow the
experimental setting of \citet{su2026online} and compare the proposed
reference-null calibrated boundary with the classical Ville boundary.

The key idea of LLM watermarking is to modify the token generation mechanism
of a language model in a way that preserves the marginal language distribution
while introducing detectable dependence on a secret random key. The detector,
which has access to the watermark key, converts the generated token sequence
into a sequence of sequential evidence variables (pivots). These pivots are
then monitored by an e-process. Under the null hypothesis of no watermark,
the pivot sequence satisfies the required calibration property, whereas
watermarked text produces systematically smaller pivots and consequently
larger e-process values.

\subsection{Language model and text generation}

We use the pretrained
\texttt{OPT-1.3B} language model without additional fine-tuning. Prompts are
sampled from the validation split of the Colossal Clean Crawled Corpus (C4).
After tokenization with the OPT tokenizer, we retain documents containing
sufficiently many tokens and use the first 50 tokens as prompts. For each
prompt, the model generates a continuation of exactly 700 tokens.

We consider two temperature settings,
$
T_{\mathrm{temp}}\in\{0.5,1\},
$
where a lower temperature corresponds to a more concentrated next-token
distribution and therefore represents a relatively weaker watermark signal,
while a higher temperature produces a more dispersed distribution and
generally makes watermark detection easier.

For each temperature, we randomly select 500 independent prompts. For every
prompt, we generate three types of text sequences:

\begin{enumerate}
    \item human-written C4 continuations following the prompt;
    \item unwatermarked OPT continuations generated by ordinary categorical
    sampling;
    \item Gumbel-max watermarked OPT continuations generated using the same
    language model and prompt.
\end{enumerate}

The first two categories serve as empirical negative controls for evaluating
false-positive behavior, while the third category is used for evaluating
watermark detection power and detection delay.

\subsection{Gumbel-max watermark and pivot construction}

Let
$
P_t=(P_{t,w})_{w\in\mathcal W}
$
denote the temperature-adjusted next-token probability vector produced by the
language model at generation step $t$, where $\mathcal W$ is the vocabulary.
For unwatermarked generation, the next token is sampled according to
\(W_t\mid P_t\sim \mathrm{Categorical}(P_t).\)
For watermarked generation, we employ the Gumbel-max watermark. At each step,
the secret watermark key generates a pseudo-random vector
$
U_t=(U_{t,w})_{w\in\mathcal W},
$
whose coordinates behave as independent uniform random variables on
$(0,1)$. The generated token is selected by
\[
W_t
=
\arg\max_{w\in\mathcal W}
\left\{
\log P_{t,w}-\log(-\log U_{t,w})
\right\}.
\]
This construction preserves the marginal token distribution:
$
P(W_t=w\mid P_t)=P_{t,w},
$
but introduces dependence between the selected token and the pseudo-random
watermark variables. Hence, the watermark does not alter the language model
distribution itself but leaves a detectable sequential signature.

Given the observed token $W_t$ and the watermark key, the detector
reconstructs the corresponding pseudo-random value
$
Y_t=U_{t,W_t},
$
and defines the one-sided pivot
$
p_t=1-Y_t .
$
Under the unwatermarked null hypothesis, the selected token is independent of
the watermark randomness, leading to
\(p_t\stackrel{\mathrm{iid}}{\sim}\mathrm{Unif}(0,1).\)
Under the Gumbel-max watermark alternative, the maximization mechanism favors
tokens associated with larger $U_{t,w}$ values. Therefore, $Y_t$ tends to be
large and $p_t$ becomes concentrated near zero.
All detection procedures are subsequently applied only to the resulting pivot
sequence
$
p_1,\ldots,p_{700}.
$

\subsection{E-process detectors}

We consider the three e-process detectors used in
\citet{su2026online}:
\(M_t^{\mathrm{WA}}, M_t^{\mathrm{OG}}, M_t^{\mathrm{AVG}}.\)
The weight-adaptive (WA) e-process uses an adaptive betting function based on
the past pivot observations. The online Grenander (OG) e-process estimates a
decreasing alternative density using the past pivots. Their average detector
is defined at the process level as
\(M_t^{\mathrm{AVG}} = \frac12M_t^{\mathrm{WA}} + \frac12M_t^{\mathrm{OG}}.\)

Importantly, the calibrated and Ville versions of each detector are evaluated
on exactly the same e-process path. Thus, the comparison isolates the effect
of the rejection boundary rather than differences in the underlying detector.

\subsection{Ville boundary and reference-null calibrated boundary}

We set the significance level to
$
\alpha=0.05.
$
The classical Ville procedure rejects whenever
$
M_t\geq \frac1\alpha=20 .
$
This threshold provides anytime-valid type-I error control but does not exploit
the finite-horizon distribution of the maximum e-process value.
To construct the reference-null calibrated boundary, we independently generate
a reference bank containing $
B=4999
$
null pivot paths, where each path consists of independent
$\mathrm{Unif}(0,1)$ observations of length 700. Each reference path is passed
through the same e-process construction, producing $
Z_b=\max_{1\leq t\leq700}M_{b,t}.
$

Let
\(Z_{(1)}\leq \cdots\leq Z_{(B)}\)
denote the ordered null maxima. The calibrated threshold is defined as
$
c_{\alpha,T}=Z_{(k_\alpha)},
$
where
$
k_\alpha=\lceil(1-\alpha)(B+1)\rceil .
$

Because the test path and reference paths are exchangeable under the null,
this construction controls the probability of at least one false rejection over
the complete monitoring horizon.

\subsection{Evaluation metrics}

We evaluate three aspects of sequential watermark detection.

First, we examine empirical sequential type-I error on human-written and
unwatermarked OPT-generated texts. For detector $m$ and boundary choice
$b\in\{C,V\}$, define
\[
\widehat{\mathrm{FPR}}_{m,b}(t)
=
\frac1{N_0}
\sum_{i=1}^{N_0}
\mathbf 1\{D_{i,b}^{(m)}\leq t\}.
\]
The main validity metric is
$\widehat{\mathrm{FPR}}_{m,b}(700)$, representing the probability of at least
one false alarm during the complete monitoring period.
Second, for watermarked texts, we evaluate detection power
\[
\widehat{\mathrm{Power}}_{m,b}(t)
=
\frac1{N_1}
\sum_{i=1}^{N_1}
\mathbf 1\{D_{i,b}^{(m)}\leq t\}.
\]
We report power curves over the whole continuation length and compare the
calibrated and Ville procedures under the paired experimental design.
Third, we evaluate detection delay. Since unsuccessful detections are
right-censored, we define
\(\widetilde D_{i,b}^{(m)} = \min\{D_{i,b}^{(m)},701\},\)
and compute the restricted mean detection delay
\(\mathrm{RMDD}_{m,b} = \frac1{N_1} \sum_{i=1}^{N_1} \widetilde D_{i,b}^{(m)} .\)
The delay gain from calibration is measured by
$
\mathrm{DG}_{m}
=
\mathrm{RMDD}_{m,V}
-
\mathrm{RMDD}_{m,C}.
$
A positive value indicates that reference-null calibration achieves earlier
watermark detection.

\subsection{Experimental results}

All results below use 500 independently sampled test sequences in each
temperature--category cell.  We write C for the reference-null calibrated
boundary and V for the Ville boundary.  A crossing is recorded when the
e-process is strictly larger than its boundary; using the same convention for
the reference bank avoids ambiguity when the pseudo-random pivots lie on a
finite computer grid.  Proportion intervals are 95\% Wilson score intervals,
and uncertainty for paired differences and delay gains is obtained from 2000
paired percentile-bootstrap resamples.  In particular, the Wilson interval is
the inversion of the binomial score test, rather than the less stable Wald
interval $\widehat p\pm1.96\sqrt{\widehat p(1-\widehat p)/n}$; this distinction
is useful here because several observed proportions are close to zero or one.

\begin{figure}[H]
    \centering
    \includegraphics[width=0.94\textwidth]{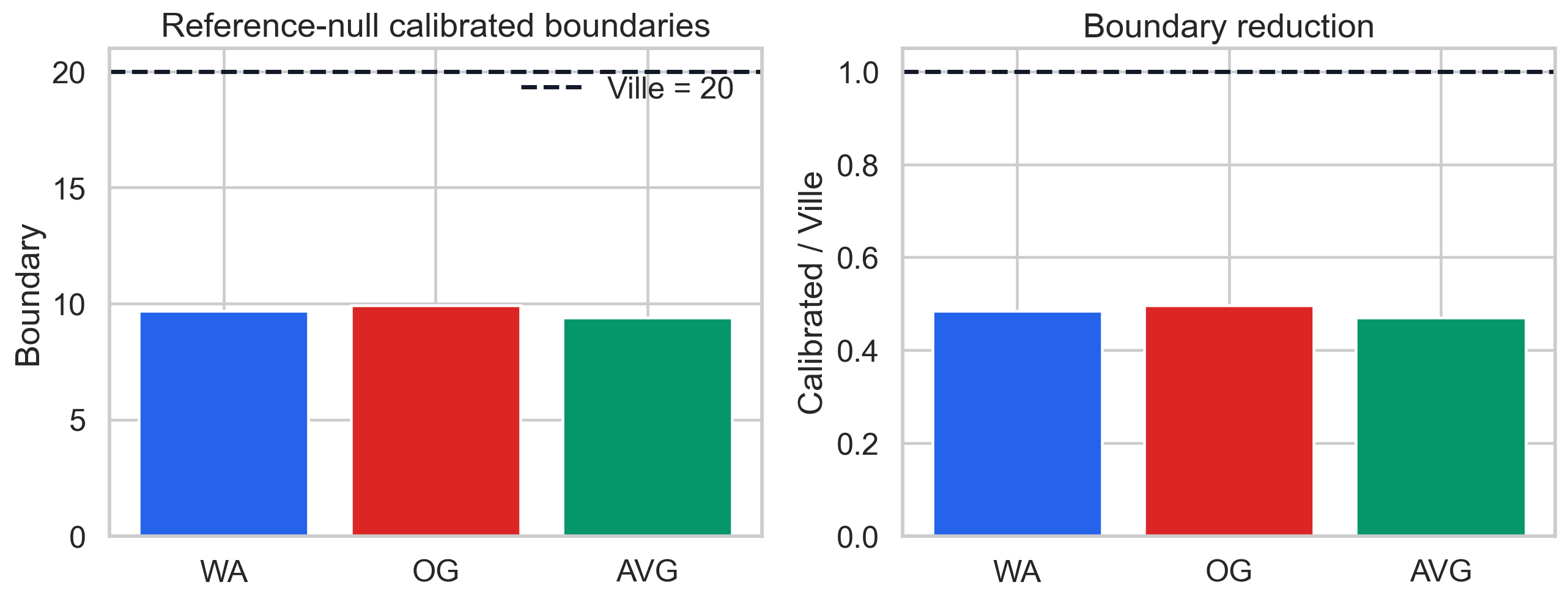}
    \caption{Calibrated boundaries and their ratios to the Ville boundary.
    Calibration reduces the required evidence by approximately one half for
    all three e-processes.}
    \label{fig:wm-boundaries}
\end{figure}

\paragraph{Calibrated boundaries.}
Table~\ref{tab:wm-boundaries} and Figure~\ref{fig:wm-boundaries} show the
detector-specific thresholds obtained from the independent bank of 4999 null
paths.  The calibrated boundary is the 4750th order statistic for every
detector.  Under strict crossing, exactly 249 of the 4999 reference maxima
exceed the selected threshold, an empirical tail probability of 0.04981.

\begin{table}[H]
    \centering
    \caption{Reference-null calibrated boundaries at $\alpha=0.05$ and
    horizon 700.  Boundary reduction is relative to the Ville value 20.}
    \label{tab:wm-boundaries}
    \begin{tabular}{lrrr}
        \toprule
        Detector & C boundary & C/V ratio & Reduction \\
        \midrule
        WA  & 9.703 & 0.485 & 51.5\%  \\
        OG  & 9.929 & 0.496 & 50.4\% \\
        AVG & 9.392 & 0.470 & 53.0\% \\
        \bottomrule
    \end{tabular}
\end{table}

The reduction is substantial and also quite stable across the three
constructions: C is between 46.96\% and 49.65\% of V.  AVG receives the
smallest absolute threshold, although this does not by itself imply that AVG
is the strongest detector because the three processes have different path
distributions.  Since the reference pivots and horizon, but not the language
model temperature, determine these thresholds, the same C values are used at
both temperatures.  This separation is important: the calibration changes
only the rejection rule and does not retrain, tune, or otherwise modify any
e-process.

\paragraph{Sequential type-I error and pivot diagnostics.}
Table~\ref{tab:wm-fpr} reports false-positive rates at the end of the 700-token
monitoring horizon.  ``Human'' denotes the untouched C4 continuation following
the prompt.  ``OPT, no WM'' denotes a continuation sampled from OPT-1.3B by
ordinary categorical sampling; the secret pseudo-random function is evaluated
only afterward to form pivots and therefore cannot influence those tokens.
These two sources provide complementary negative controls for natural and
model-generated text.

\begin{table}[H]
    \centering
    \caption{Sequential false-positive rate in percent at token 700.  Brackets
    contain 95\% Wilson score intervals; each row uses $n=500$ sequences.}
    \label{tab:wm-fpr}
    \scriptsize
    \setlength{\tabcolsep}{4.2pt}
    \begin{tabular}{clcll}
        \toprule
        Temp. & Negative control & Detector & C: FPR [95\% CI] & V: FPR [95\% CI] \\
        \midrule
        0.5 & Human      & WA  & $5.0\ [3.4,7.3]$ & $2.6\ [1.5,4.4]$ \\
        0.5 & Human      & OG  & $4.4\ [2.9,6.6]$ & $1.8\ [0.9,3.4]$ \\
        0.5 & Human      & AVG & $4.4\ [2.9,6.6]$ & $2.0\ [1.1,3.6]$ \\
        0.5 & OPT, no WM & WA  & $6.0\ [4.2,8.4]$ & $2.8\ [1.7,4.6]$ \\
        0.5 & OPT, no WM & OG  & $5.6\ [3.9,8.0]$ & $2.0\ [1.1,3.6]$ \\
        0.5 & OPT, no WM & AVG & $5.6\ [3.9,8.0]$ & $1.8\ [0.9,3.4]$ \\
        \midrule
        1.0 & Human      & WA  & $4.8\ [3.2,7.0]$ & $3.2\ [2.0,5.1]$ \\
        1.0 & Human      & OG  & $5.2\ [3.6,7.5]$ & $2.4\ [1.4,4.1]$ \\
        1.0 & Human      & AVG & $5.6\ [3.9,8.0]$ & $2.8\ [1.7,4.6]$ \\
        1.0 & OPT, no WM & WA  & $3.0\ [1.8,4.9]$ & $1.0\ [0.4,2.3]$ \\
        1.0 & OPT, no WM & OG  & $5.2\ [3.6,7.5]$ & $1.4\ [0.7,2.9]$ \\
        1.0 & OPT, no WM & AVG & $4.2\ [2.8,6.3]$ & $1.2\ [0.6,2.6]$ \\
        \bottomrule
    \end{tabular}
\end{table}

\begin{figure}[H]
    \centering
    \includegraphics[width=0.97\textwidth]{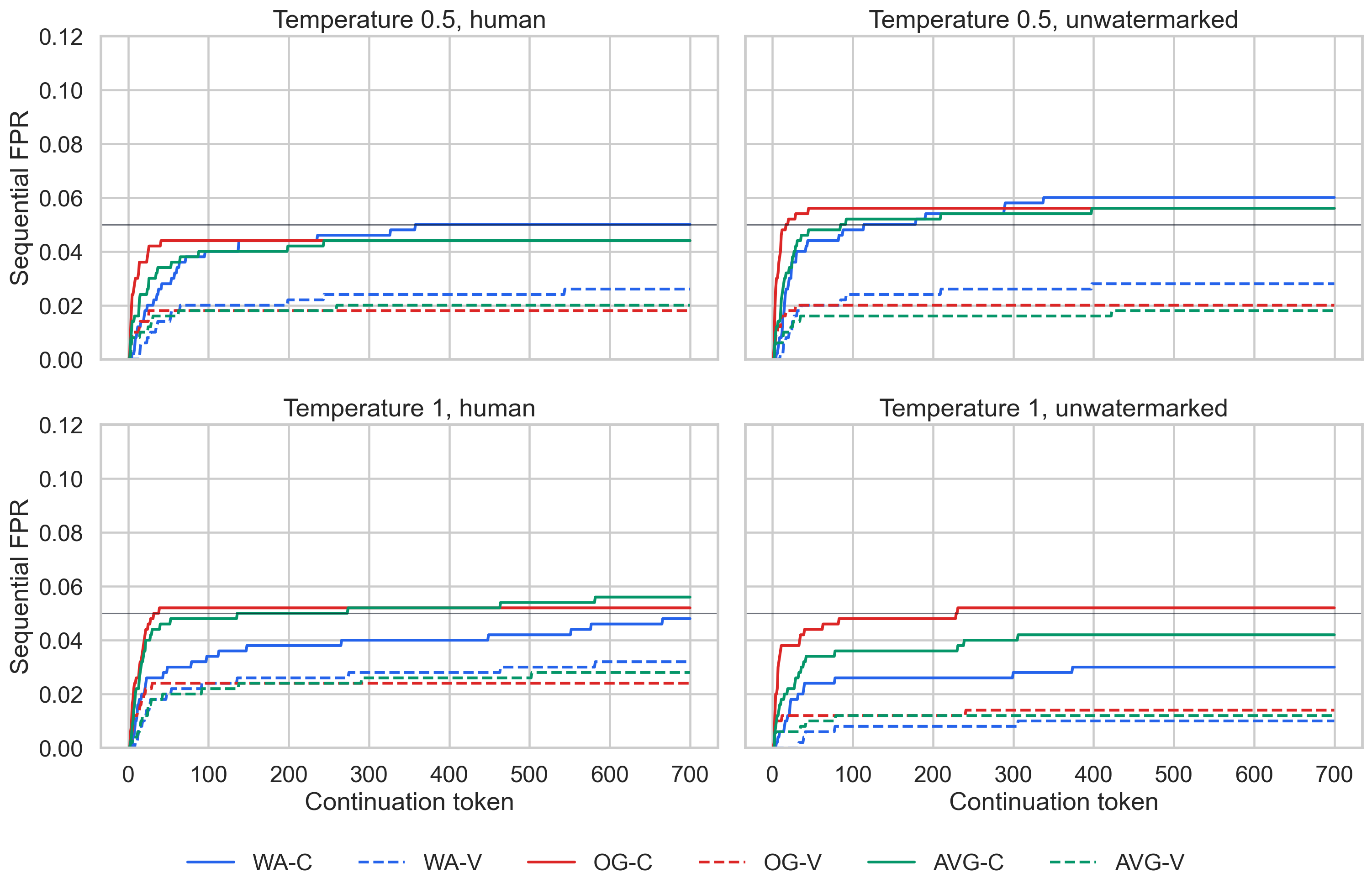}
    \caption{Cumulative sequential false-positive rates.  Solid lines use C,
    dashed lines use V, and the horizontal line is the nominal 5\% level.  A
    curve records whether each sequence has crossed by the displayed token, so
    it is necessarily nondecreasing.}
    \label{fig:wm-fpr}
\end{figure}

Across the 12 detector--control combinations, C produces endpoint FPRs between
3.0\% and 6.0\%.  Eleven of the 12 Wilson intervals contain the nominal 5\%
level; the sole exception is WA on unwatermarked OPT text at temperature 1,
where the estimate is conservative (3.0\%, interval $[1.8\%,4.9\%]$).  There
is consequently no empirical indication that the roughly twofold boundary
reduction causes systematic inflation beyond the target level.  By contrast,
all V point estimates are between 1.0\% and 3.2\%.  This is valid but visibly
conservative, and is precisely the unused type-I error budget that calibration
converts into earlier rejection.

\begin{figure}[H]
    \centering
    \includegraphics[width=0.96\textwidth]{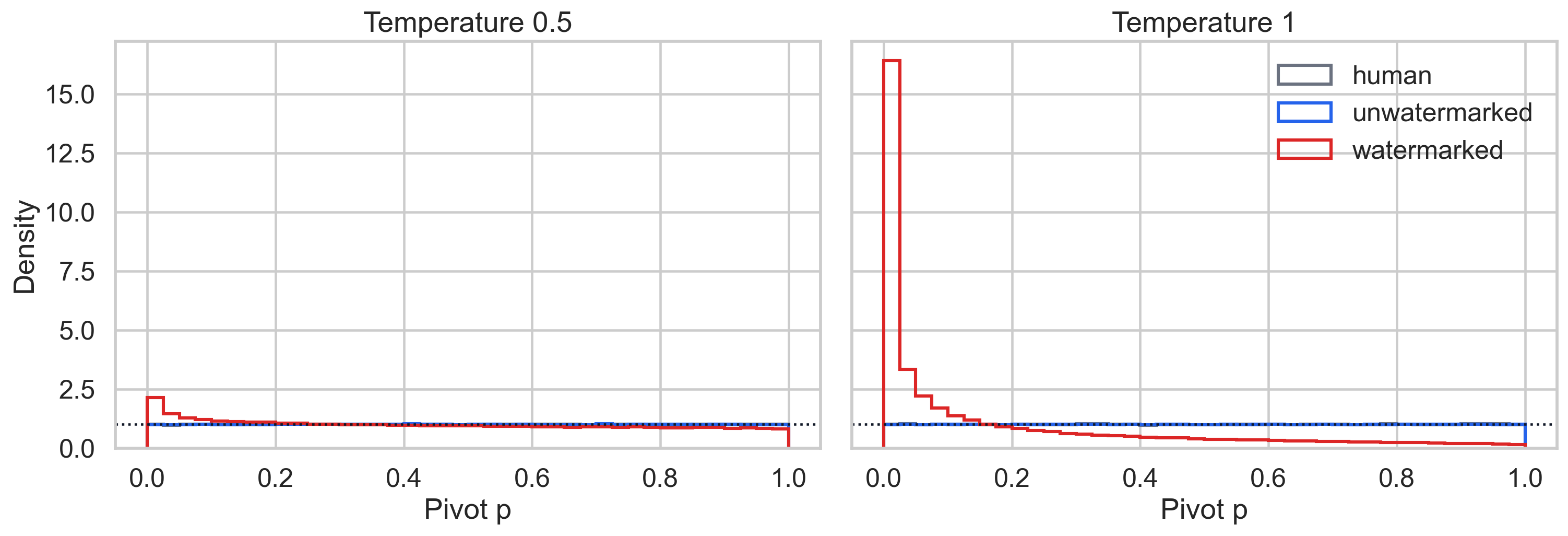}
    \caption{Empirical pivot densities for human, unwatermarked OPT, and
    watermarked OPT text.  The dotted line is the $\Unif(0,1)$ density.  The
    negative controls remain nearly uniform, whereas watermarking shifts mass
    toward zero, especially at temperature 1.}
    \label{fig:wm-pivots}
\end{figure}

The pathwise curves in Figure~\ref{fig:wm-fpr} also show that most false alarms
occur early, after which the estimates stabilize.  The similarity of the human
and unwatermarked OPT controls argues against a model-specific artifact.  This
conclusion is reinforced by the pivot diagnostics in
Figure~\ref{fig:wm-pivots}: over 350,000 pivots per cell, the null means range
only from 0.49964 to 0.50037, and Kolmogorov--Smirnov tests against the uniform
distribution give $p$-values from 0.311 to 0.583.  These diagnostics are not a
substitute for the exchangeability argument, but they verify that the
implemented key reconstruction and null data pipeline behave as intended.

\begin{figure}[H]
    \centering
    \includegraphics[width=0.97\textwidth]{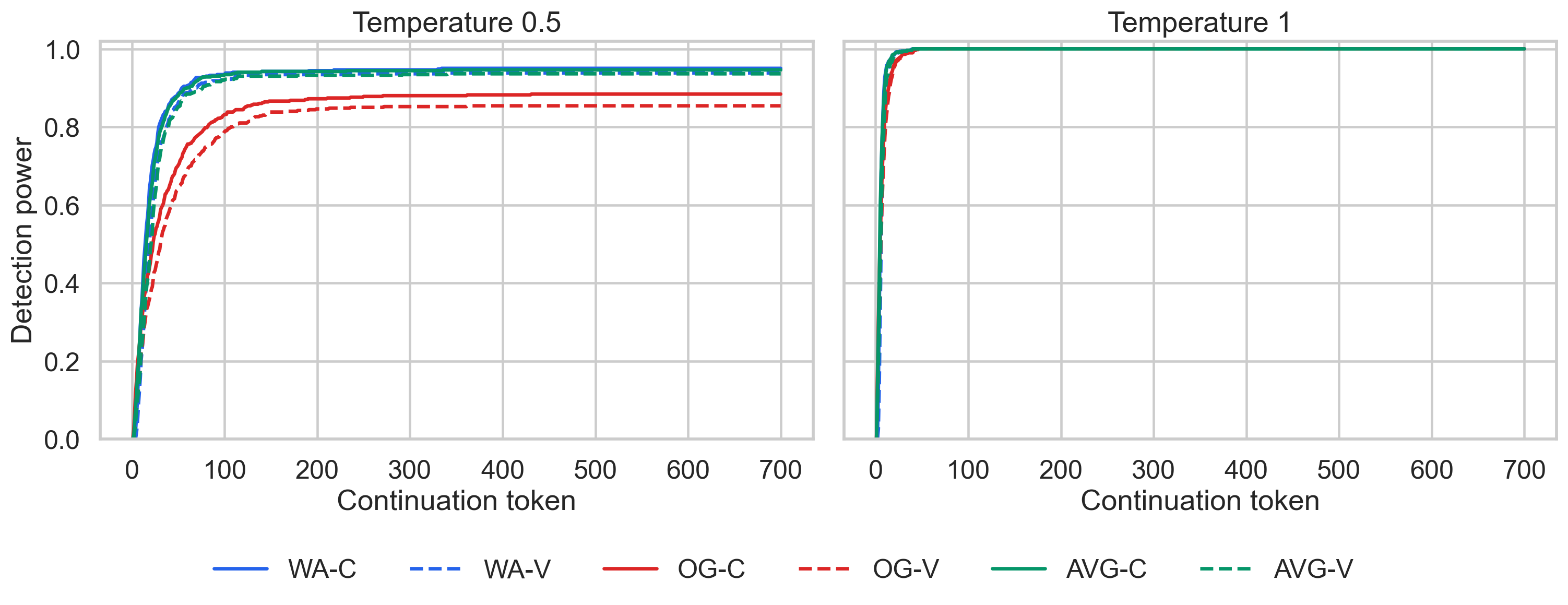}
    \caption{Cumulative detection power under watermarking.  Solid and dashed
    lines denote C and V, respectively.  Calibration weakly dominates Ville
    pathwise because it applies a lower threshold to the same e-process.}
    \label{fig:wm-power}
\end{figure}

\paragraph{Detection power.}
Table~\ref{tab:wm-power} summarizes early and full-horizon power, and
Figure~\ref{fig:wm-power} gives the complete power curves.  Full-horizon
intervals are Wilson intervals for the individual power estimates; intervals
for C--V are paired bootstrap intervals and therefore exploit the fact that
both boundaries are applied to the same e-process paths.

\begin{table}[H]
    \centering
    \caption{Watermark detection power in percent at tokens 100 and 700.
    Brackets at token 700 are 95\% intervals.  The last column is the paired
    percentage-point gain of C over V at token 700.}
    \label{tab:wm-power}
    \scriptsize
    \resizebox{\textwidth}{!}{%
    \begin{tabular}{clrrrrr}
        \toprule
        Temp. & Detector & C at 100 & V at 100 & C at 700 [95\% CI] & V at 700 [95\% CI] & C--V [95\% CI] \\
        \midrule
        0.5 & WA  & 93.4 & 92.2 & $95.0\ [92.7,96.6]$ & $93.8\ [91.3,95.6]$ & $1.2\ [0.4,2.2]$ \\
        0.5 & OG  & 83.2 & 78.8 & $88.4\ [85.3,90.9]$ & $85.4\ [82.0,88.2]$ & $3.0\ [1.6,4.6]$ \\
        0.5 & AVG & 93.4 & 92.0 & $94.6\ [92.3,96.3]$ & $93.6\ [91.1,95.4]$ & $1.0\ [0.2,2.0]$ \\
        \midrule
        1.0 & WA  & 100.0 & 100.0 & $100.0\ [99.2,100.0]$ & $100.0\ [99.2,100.0]$ & $0.0\ [0.0,0.0]$ \\
        1.0 & OG  & 100.0 & 100.0 & $100.0\ [99.2,100.0]$ & $100.0\ [99.2,100.0]$ & $0.0\ [0.0,0.0]$ \\
        1.0 & AVG & 100.0 & 100.0 & $100.0\ [99.2,100.0]$ & $100.0\ [99.2,100.0]$ & $0.0\ [0.0,0.0]$ \\
        \bottomrule
    \end{tabular}}
\end{table}

At temperature 0.5, where the watermark signal is weaker, the benefit of C is
visible throughout the monitoring period.  At token 50 the gains are 2.4,
6.0, and 2.8 percentage points for WA, OG, and AVG, with paired 95\% intervals
$[1.2,3.8]$, $[4.0,8.2]$, and $[1.6,4.2]$, respectively.  The advantage
narrows as the easier paths are detected by both procedures, but remains
positive at token 700 for all three detectors.  OG has lower absolute power
than WA and AVG, but obtains the largest calibration gain: 4.4 percentage
points at token 100 and 3.0 points at token 700.  Thus, a detector need not be
the strongest in absolute terms to benefit most from a less conservative
boundary.

At temperature 1, every watermarked sequence is detected by every method and
boundary by token 50.  This ceiling explains the zero observed power
difference; it does not mean that the boundary choice is irrelevant, because
the crossing times can still differ.  The pivot means clarify the temperature
effect: the watermarked mean decreases from 0.45420 at temperature 0.5 to
0.17649 at temperature 1, while the null means remain at 0.5.  A more diffuse
next-token distribution permits the keyed Gumbel variables to exert greater
influence on the selected token, producing much smaller pivots and much faster
evidence accumulation.



\paragraph{Detection delay.}
Power at a fixed horizon conceals gains on paths that both methods eventually
detect.  Table~\ref{tab:wm-delay} therefore reports the restricted mean
detection delay with nondetections coded as 701, as specified above.

\begin{table}[H]
    \centering
    \caption{Restricted mean detection delay (tokens) and the gain from
    calibration.  Confidence intervals are based on 2000 paired bootstrap
    resamples.}
    \label{tab:wm-delay}
    \scriptsize
    \resizebox{\textwidth}{!}{%
    \begin{tabular}{clrrrr}
        \toprule
        Temp. & Detector & RMDD(C) & RMDD(V) & V--C gain [95\% CI] & Relative reduction [95\% CI] \\
        \midrule
        0.5 & WA  & 55.31  & 66.49  & $11.17\ [6.22,17.48]$  & $16.8\%\ [9.7,25.0]$ \\
        0.5 & OG  & 111.20 & 135.08 & $23.88\ [15.17,34.17]$ & $17.7\%\ [11.8,24.6]$ \\
        0.5 & AVG & 57.90  & 68.71  & $10.81\ [5.53,17.26]$  & $15.7\%\ [8.4,23.6]$ \\
        \midrule
        1.0 & WA  & 6.32 & 7.50 & $1.19\ [1.08,1.30]$ & $15.8\%\ [14.4,17.2]$ \\
        1.0 & OG  & 7.22 & 8.09 & $0.87\ [0.69,1.08]$ & $10.7\%\ [8.7,13.2]$ \\
        1.0 & AVG & 5.97 & 7.04 & $1.06\ [0.91,1.22]$ & $15.1\%\ [13.2,17.1]$ \\
        \bottomrule
    \end{tabular}}
\end{table}

\begin{figure}[H]
    \centering
    \includegraphics[width=0.94\textwidth]{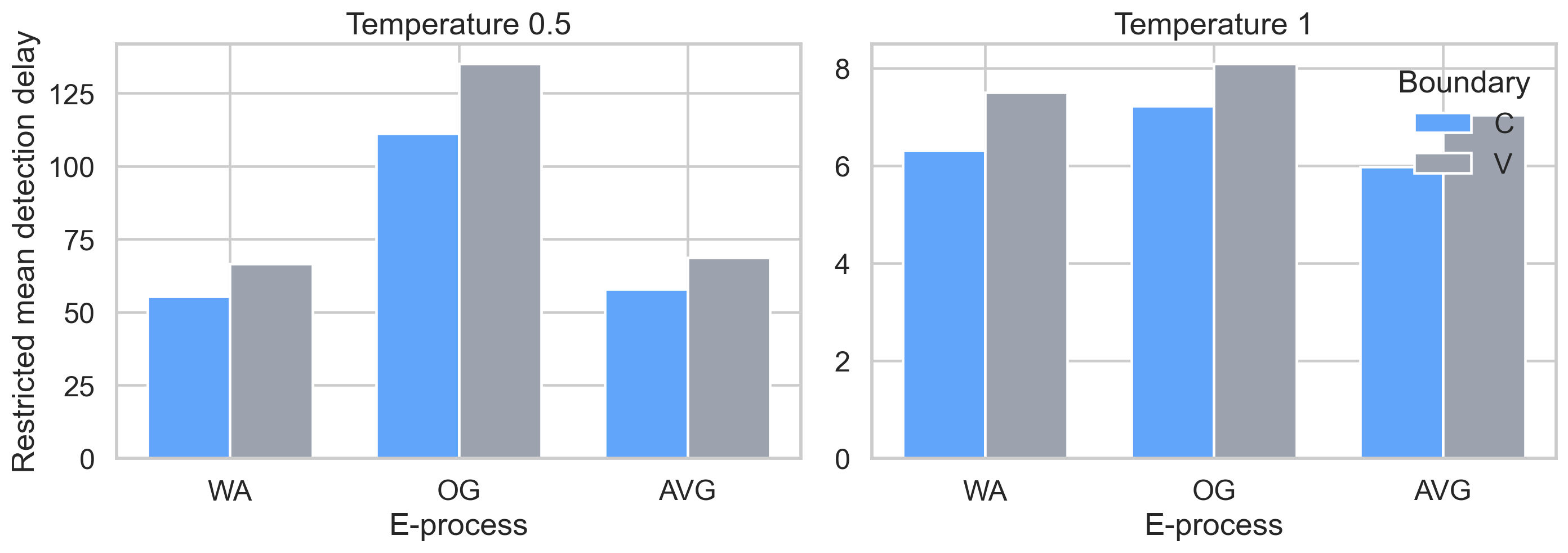}
    \caption{Restricted mean detection delay under C and V.  The two panels
    use different vertical scales because the temperature-1 watermark is
    detected within only a few tokens.}
    \label{fig:wm-delay}
\end{figure}

All six paired delay gains are positive and all six bootstrap intervals exclude
zero.  At temperature 0.5, calibration saves approximately 11 tokens for WA
and AVG and 24 tokens for OG; relative RMDD reductions lie between 15.7\% and
17.7\%.  At temperature 1 the absolute savings are necessarily smaller because
both methods are already extremely fast, yet C still reduces RMDD by
10.7\%--15.8\%.  AVG has the smallest RMDD at both temperatures, whereas OG
has the largest absolute calibration gain in the difficult temperature-0.5
condition.

Taken together, the results support the intended modular interpretation of
reference-null calibration.  The calibrated boundaries spend nearly all of
the available 5\% sequential error budget on two distinct negative controls,
while the Ville boundaries leave a substantial portion unused.  Applying the
smaller thresholds to unchanged WA, OG, and AVG paths yields uniformly earlier
pathwise crossings, statistically resolved delay reductions, and positive
power gains whenever the problem is not already at the ceiling.  The numerical
conclusions are specific to OPT-1.3B, C4 prompts, the two temperatures, and the
700-token horizon; extending them to other models, watermark keys, decoding
schemes, or horizons requires a new independent reference bank matched to the
prospective detector and monitoring protocol.

\section{Conclusion}
\label{sec:conclusion}

In this paper, we developed a reference-null calibration framework for
constructing sharper rejection thresholds for e-values and finite-horizon
e-processes. Rather than relying on the universal Markov or Ville boundary,
the proposed approach uses independent reference-null realizations to
calibrate the relevant null distribution directly. The resulting procedure
leaves the underlying e-value or e-process unchanged and separates evidence
construction from boundary calibration. Under exchangeability between the
test statistic and the reference-null statistics, the calibrated rule retains
finite-sample type-I error control while adapting the rejection threshold to
the actual null behavior of the procedure.

We studied the statistical consequences of this calibration principle in
detail for conformal martingales. For adaptive histogram betting, we used
Krichevsky--Trofimov smoothing and considered both the unrestarted process and
a restart-mixture construction for unknown change points. The finite-sample
analysis shows explicitly how sequential detection performance depends on
the strength of the rank signal, histogram discretization, online learning,
and restart aggregation. In particular, the results quantify how replacing
the Ville boundary by a smaller reference-null calibrated boundary reduces
the evidence required for detection and can improve both detection power and
detection delay, without modifying the underlying betting strategy.

The numerical experiments support these theoretical conclusions. Across the
distribution-shift settings considered in the simulations, reference-null
calibration improves sequential detection when the universal Ville boundary
is conservative. The LLM watermark detection experiment further illustrates
the modular nature of the approach: the same WA, OG, and AVG e-processes are
used under both procedures, and the improvement is obtained solely by
recalibrating their rejection thresholds. This makes reference-null
calibration particularly attractive when a reliable e-process already exists
but its universal rejection boundary leaves a substantial portion of the
available type-I error budget unused.

Several limitations also suggest directions for future work. The present
framework requires a reference sample that accurately represents the null
distribution relevant to the prospective test, and the finite-horizon
calibrated boundary depends on the monitoring horizon and other design
choices fixed before calibration. Extending the methodology to composite
null hypotheses without a pivotal reduction, approximate rather than exact
reference-null exchangeability, and settings in which calibration data must
be reused across multiple monitoring problems would broaden its applicability.
It would also be useful to study more adaptive boundary shapes and restart
schemes while retaining finite-sample validity. More generally, the results
suggest that external null information can be used not only to construct
e-processes, but also to improve the efficiency of already valid sequential
procedures through a separate and reusable calibration layer.

\clearpage

\bibliographystyle{plainnat}
\bibliography{references}

@article{vovk2021evalues,
  author  = {Vovk, Vladimir and Wang, Ruodu},
  title   = {E-values: Calibration, combination, and applications},
  journal = {The Annals of Statistics},
  volume  = {49},
  number  = {3},
  pages   = {1736--1754},
  year    = {2021},
  doi     = {10.1214/20-AOS2020}
}

@article{shafer2021testing,
  author  = {Shafer, Glenn},
  title   = {Testing by betting: A strategy for statistical and scientific communication},
  journal = {Journal of the Royal Statistical Society: Series A (Statistics in Society)},
  volume  = {184},
  number  = {2},
  pages   = {407--431},
  year    = {2021},
  doi     = {10.1111/rssa.12647}
}

@article{grunwald2024safe,
  author  = {Gr{\"u}nwald, Peter and de Heide, Rianne and Koolen, Wouter M.},
  title   = {Safe testing},
  journal = {Journal of the Royal Statistical Society Series B: Statistical Methodology},
  volume  = {86},
  number  = {5},
  pages   = {1091--1128},
  year    = {2024},
  doi     = {10.1093/jrsssb/qkae011}
}

@article{ramdas2023game,
  author  = {Ramdas, Aaditya and Gr{\"u}nwald, Peter and Vovk, Vladimir
             and Shafer, Glenn},
  title   = {Game-Theoretic Statistics and Safe Anytime-Valid Inference},
  journal = {Statistical Science},
  volume  = {38},
  number  = {4},
  pages   = {576--601},
  year    = {2023},
  doi     = {10.1214/23-STS894}
}

@article{blierwong2024improved,
  author  = {Blier-Wong, Christopher and Wang, Ruodu},
  title   = {Improved thresholds for e-values},
  journal = {arXiv preprint arXiv:2408.11307},
  year    = {2024},
  note    = {Forthcoming in The Annals of Statistics}
}

@article{fischer2026improving,
  author  = {Fischer, Lasse and Ramdas, Aaditya},
  title   = {Improving Wald's (Approximate) Sequential Probability Ratio
             Test by Avoiding Overshoot},
  journal = {IEEE Transactions on Information Theory},
  volume  = {72},
  number  = {4},
  pages   = {2457--2471},
  year    = {2026}
}

@article{delapena2026exact,
  author  = {de la Pe{\~n}a, Victor H. and Klass, Michael J.},
  title   = {The Exact Ville Identity: From the Absorbing Case to the
             General Law with an Application to E-Values},
  journal = {arXiv preprint arXiv:2607.04620},
  year    = {2026}
}

@inproceedings{vovk2003testing,
  author    = {Vovk, Vladimir and Nouretdinov, Ilia and Gammerman, Alex},
  title     = {Testing Exchangeability On-Line},
  booktitle = {Proceedings of the 20th International Conference on Machine Learning},
  pages     = {768--775},
  year      = {2003},
  publisher = {AAAI Press}
}

@book{vovk2005algorithmic,
  author    = {Vovk, Vladimir and Gammerman, Alexander and Shafer, Glenn},
  title     = {Algorithmic Learning in a Random World},
  publisher = {Springer},
  address   = {New York},
  year      = {2005}
}

@inproceedings{fedorova2012plugin,
  author    = {Fedorova, Valentina and Gammerman, Alex and Nouretdinov, Ilia
               and Vovk, Vladimir},
  title     = {Plug-in Martingales for Testing Exchangeability On-Line},
  booktitle = {Proceedings of the 29th International Conference on Machine Learning},
  pages     = {1639--1646},
  year      = {2012},
  publisher = {Omnipress}
}

@inproceedings{eliades2020histogram,
  author    = {Eliades, Charalambos and Papadopoulos, Harris},
  title     = {A Histogram Based Betting Function for Conformal Martingales},
  booktitle = {Proceedings of the Ninth Symposium on Conformal and
               Probabilistic Prediction and Applications},
  series    = {Proceedings of Machine Learning Research},
  volume    = {128},
  pages     = {100--113},
  year      = {2020},
  publisher = {PMLR}
}

@article{krichevsky1981performance,
  author  = {Krichevsky, Raphail E. and Trofimov, Victor K.},
  title   = {The Performance of Universal Encoding},
  journal = {IEEE Transactions on Information Theory},
  volume  = {27},
  number  = {2},
  pages   = {199--207},
  year    = {1981},
  doi     = {10.1109/TIT.1981.1056331}
}

@article{su2026online,
  author  = {Su, Weijie and Wang, Ruodu and Zhao, Zinan},
  title   = {Online {LLM} Watermark Detection via E-Processes},
  journal = {arXiv preprint arXiv:2602.14286},
  year    = {2026}
}

@article{lei2018distribution,
  author  = {Lei, Jing and G'Sell, Max and Rinaldo, Alessandro
             and Tibshirani, Ryan J. and Wasserman, Larry},
  title   = {Distribution-Free Predictive Inference for Regression},
  journal = {Journal of the American Statistical Association},
  volume  = {113},
  number  = {523},
  pages   = {1094--1111},
  year    = {2018},
  doi     = {10.1080/01621459.2017.1307116}
}

@article{barber2021predictive,
  author  = {Barber, Rina Foygel and Cand{\`e}s, Emmanuel J.
             and Ramdas, Aaditya and Tibshirani, Ryan J.},
  title   = {Predictive Inference with the Jackknife+},
  journal = {The Annals of Statistics},
  volume  = {49},
  number  = {1},
  pages   = {486--507},
  year    = {2021},
  doi     = {10.1214/20-AOS1965}
}

@article{saha2026nonpartitioned,
  title={Non-partitioned e-detectors for nonparametric sequential change detection},
  author={Saha, Aytijhya and Ramdas, Aaditya},
  journal={arXiv preprint arXiv:2607.28322},
  year={2026},
  month={July},
  archivePrefix={arXiv},
  eprint={2607.28322},
  primaryClass={stat.ME}
}

@inproceedings{shaer2026cctm,
  title     = {Testing For Distribution Shifts with Conditional Conformal Test Martingales},
  author    = {Shaer, Shalev and Bar, Yarin and Prinster, Drew and Romano, Yaniv},
  booktitle = {Proceedings of the 43rd International Conference on Machine Learning},
  year      = {2026},
  note      = {arXiv:2602.13848}
}

@article{vovk2019nonparametric,
  title   = {Nonparametric Predictive Distributions Based on Conformal Prediction},
  author  = {Vovk, Vladimir and Shen, Jieli and Manokhin, Valery and Xie, Min-ge},
  journal = {Machine Learning},
  volume  = {108},
  pages   = {445--474},
  year    = {2019},
  doi     = {10.1007/s10994-018-5755-8}
}

@article{mann1947test,
  title={On a Test of Whether One of Two Random Variables Is Stochastically Larger than the Other},
  author={Mann, H. B. and Whitney, D. R.},
  journal={The Annals of Mathematical Statistics},
  volume={18},
  number={1},
  pages={50--60},
  year={1947},
  doi={10.1214/aoms/1177730491}
}

@article{hoeffding1948class,
  title={A Class of Statistics with Asymptotically Normal Distribution},
  author={Hoeffding, Wassily},
  journal={The Annals of Mathematical Statistics},
  volume={19},
  number={3},
  pages={293--325},
  year={1948},
  doi={10.1214/aoms/1177730196}
}

@book{hajek1999theory,
  title={Theory of Rank Tests},
  author={H{\'a}jek, Jaroslav and {\v S}id{\'a}k, Zbyn{\v e}k and Sen, Pranab K.},
  edition={2},
  publisher={Academic Press},
  year={1999}
}

@inproceedings{eliades2022betting,
  title     = {A Betting Function for Addressing Concept Drift with Conformal Martingales},
  author    = {Eliades, Charalambos and Papadopoulos, Harris},
  booktitle = {Proceedings of the Eleventh Symposium on Conformal and Probabilistic Prediction with Applications},
  series    = {Proceedings of Machine Learning Research},
  volume    = {179},
  pages     = {219--238},
  year      = {2022},
  publisher = {PMLR}
}

@article{farran2026model,
  title   = {When Your Model Stops Working: Anytime-Valid Calibration Monitoring},
  author  = {Farran, Tristan},
  journal = {arXiv preprint arXiv:2603.13156},
  year    = {2026}
}

@book{wald1945sequential,
  title={Sequential Analysis},
  author={Wald, Abraham},
  year={1945},
  publisher={John Wiley \& Sons},
  address={New York}
}

@article{lai1988nearly,
  title={Nearly optimal sequential tests of composite hypotheses},
  author={Lai, Tze Leung},
  journal={The Annals of Statistics},
  volume={16},
  number={3},
  pages={856--886},
  year={1988},
  publisher={Institute of Mathematical Statistics}
}

@inproceedings{kharitonov2015sequential,
  title={Sequential Testing for Early Stopping of Online Experiments},
  author={Kharitonov, Eugene and Vorobev, Artem and Ustinovskiy, Yury and Serdyukov, Pavel},
  booktitle={Proceedings of the 21st ACM SIGKDD International Conference on Knowledge Discovery and Data Mining},
  pages={1977--1986},
  year={2015},
  publisher={ACM}
}

@article{waudby2023distribution,
  title={Distribution-uniform anytime-valid sequential inference},
  author={Waudby-Smith, Ian and Kennedy, Edward H. and Ramdas, Aaditya},
  journal={arXiv preprint arXiv:2311.03343},
  year={2023}
}

@article{lindon2024anytime,
  title={Anytime-Valid Inference in Linear Models and Regression-Adjusted Experiments},
  author={Lindon, Matthew and others},
  journal={Harvard Business School Working Paper},
  number={24-060},
  year={2024}
}

\end{document}